\documentclass[11pt]{article}
\usepackage{amsmath,amssymb}
\usepackage{enumerate}
\usepackage{color}

\usepackage{xy}
\usepackage{xypic}
\usepackage{subfigure}
\usepackage{tikz-cd}

\newcommand{\R}{\mathbb{R}}

\newcommand{\Q}{\mathbb{Q}}
\newcommand{\N}{\mathbb{N}}
\newcommand{\Z}{\mathbb{Z}}
\newcommand{\bea}{\begin{eqnarray}}
\newcommand{\eea}{\end{eqnarray}}
\newcommand{\la}{\lambda}
\newcommand{\ccD}{\check D}

\def\a{\alpha}

\def\de{\delta}
\def\e{\varepsilon}

\def\s{\sigma}

\def\diam{{\rm diam}}

\def\1{\rm Id}
\def\sgn{{\rm sgn}}
\def\sup{{\rm sup}}

\def\ci{\circ}
\def\cl{{\rm cl}}
\def\pr{{\rm pr}}
\def\tdiam{{\rm tdiam}}
\newcommand{\qed}{$\hfill\blacksquare$}

\def\V{\noindent}

\def\MCS{\mathcal{MCS}}

\def\cl{{\rm cl}}

\def\CPOM{\rm {\bf POM}}

\def\CALL{\rm {\bf ALL}}

\def\dist{{\rm dist}}
\def\Corr{{\rm Corr}}
\def\Dm{{\rm Dm}}

\newcommand{\bean}{\begin{eqnarray*}}
\newcommand{\eean}{\end{eqnarray*}}

\newtheorem{Theorem}{Theorem}

\newtheorem{Definition}{Definition}

\newcommand{\ben}{\begin{enumerate}}
\newcommand{\een}{\end{enumerate}}
\newcommand{\bit}{\begin{itemize}}
\newcommand{\eit}{\end{itemize}}
\newcommand{\edoc}{\end{document}}

\usepackage{hyperref}

\usepackage{enumerate}

\title{\vspace{-1.8cm}Gromov-Hausdorff metrics and dimensions of Lorentzian length spaces}

\begin{document}

\author{Olaf M\"uller\footnote{Institut f\"ur Mathematik, Humboldt-Universit\"at zu Berlin, Unter den Linden 6, D-10099 Berlin, \texttt{Email: o.mueller@hu-berlin.de}}}

\date{\today}
\maketitle

\begin{abstract}
\V We construct analoga of Gromov-Hausdorff space for Lorentzian distances and show a Gromov precompactness result for one of them. After calculating the Dushnik-Miller dimension of Minkowski spaces (of manifold dimension larger than 2) to be countable infinity, we define a dimension for ordered sets recovering the correct manifold dimension, obtain an obstruction for existence of injective monotonous maps between Lorentzian length spaces, induce functorial pseudo-metrics on Cauchy subsets that in the spacetime case coincide with the Riemannian ones, and prove existence of anti-Lipschitz Cauchy functions with a given Cauchy zero locus, a fundamental ingredient for the Sormani-Vega null distance.
\end{abstract}

\section{Introduction}

Recently, attempts to 'synthesize' Lorentzian geometry in a similar way as in Riemannian geometry have attracted much attention. Sec. \ref{GHConv} of this article presents and compares three notions of Gromov-Hausdorff convergence for Lorentzian length spaces and shows a Gromov precompactness result for one of them, $d_{GH}^-$. In Sec. \ref{Dim}, motivated by the definition of the Hausdorff dimension for Lorentzian length spaces \cite{rMcS}, we explore different notions of dimension in the setting of (almost) Lorentzian length spaces. After a short discussion of the notion of dimension in different categories, the focus is on notions of dimension using only the order relation. Then we calculate 
$\dim_{DM}(\R^{1,n}) = \aleph_0 \  \forall n \geq 2$
for the classical Dushnik-Miller dimension $\dim_{DM}$ of an ordered set and $\R^{1,n}$ being the $(n+1)$-dimensional Minkowski spacetime. For several other previously defined order dimensions $\dim_*$ (see \cite{HBG}) we show a lower estimate $\dim_*(\R^{1,n}) \geq \aleph_0$ for all $n \geq 2$. 
To remedy this incompatibility between the manifold dimension and the known order dimensions, we present a novel notion of dimension for ordered sets that takes the expected value $n+1$ on $\R^{1,n}$ and are lower semi-continuous with respect to the Lorentzian Gromov-Hausdorff metric $d_{GH}^-$ defined before. In Sec. 3, we induce functorial pseudometrics on Cauchy subsets that in the spacetime case coincide with the Riemannian ones, prove existence of anti-Lipschitz Cauchy functions with a given Cauchy zero locus, which is a fundamental ingredient for the Sormani-Vega null distance \cite{SV}, and consider collapse phenomena in this context. 

\newpage

\section{Notions of convergence of Lorentzian length spaces}
\label{GHConv}

Let $S$ be an arbitrary set fixed once and forever in this article. For all categories in this article we infer without further mention that all objects are of cardinality $\leq \# S$. For a category ${\bf C}$, let $C^I$ denote the set of isomorphism classes of objects of ${\bf C}$. For a functor $\mathcal{F}: {\bf X} \rightarrow {\bf Y}$ between two categories, let $\mathcal{F}^I : X^I \rightarrow Y^I$ be the push-down of $\mathcal{F}$. We denote the closure of $A \subset X$ in a topological space $X$ by $\cl (A,X)$ and the set-theoretic symmetric difference of two subsets $A,B$ of a set $C$ by $A \triangle B$.

\bigskip

In synthetic Lorentzian geometry, there are two different functorial approaches. The first one uses the category $\CALL$ of (causal) almost Lorentzian pre-length spaces (\cite{mKcS}, \cite{oM:LorFunct}), which are sets $X$ with a signed Lorentzian distance function $\s$ on $X$, i.e., a function $\s: X \times X \rightarrow \R$ that is antisymmetric and satisfies $\s (x,z) \geq \s (x,y) + \s(y,z) $ whenever $\s(x,y), \s(y,z) >0$. \footnote{In the above references, $\s$ is moreover supposed to be lower semi-continuous w.r.t. a topology given or defined from $\s$, like the topology $\tau_+$ described below; here, however, we will add continuity to our assumptions explicitly if we need it. Moreover, in our context, $\s$ takes values of both signs, whereas many authors prefer a function $\tau : X \times X \rightarrow \R\cup \{ - \infty \}$ supported on $J^+$, zero on $J^+ \setminus I^+ $ and taking the value $- \infty$ on $X \times X \setminus J^+ $. The advantage of this definition is that we can read off immediately not only $I^+$ but also $J^+$ from $\tau$. However, if the causal relation $\leq$ is supposed to be related to the chronological relation by e.g. $ x \leq y \Leftrightarrow I^+(y) \subset I^+(x) \lor I^- (x) \subset I^- (y)$, both approaches are easily seen to be essentially equivalent. Moreover, we prefer to use the letter "$\s$" here to reserve the letter "$\tau$" for topologies like $\tau_+$.} On Lorentzian manifolds, there is a natural signed Lorentzian distance function $\s (p,q) := \pm \sup \{ \ell (c) | c: p \leadsto q \in J^{\pm} (p) {\ \rm  causal} \} $ where $\ell (c) \geq 0$ is the length of a Lipschitz causal curve defined in analogy to metric length, replacing $\sup $ with $\inf$. The other approach uses the category $\CPOM$ of (partially) ordered measure spaces. The two categories seem to be intimately connected (see below, see \cite{oM:LorFunct} for further details). 

Let $(X, \leq)$ be an ordered space. We can define a chronological relation $\ll$ on $X$ by $\beta$ or $\gamma$ as follows: 

Firstly, following a suggestion by Miguzzi and S\'anchez \cite{eMmS}, we define a binary relation $ \beta ( \leq )$ via 

$$(x,y) \in \beta ( \leq ) : \Leftrightarrow \big( x \leq y \land ( \exists u,v \in X:  x < u< v < y \land J^+(u) \cap J^-(v) {\rm \ not \ totally \ ordered} ) \big)$$

Secondly, we define 

$$ p (\gamma(\leq) ) q : \Leftrightarrow (\forall a>p  \exists  a>b>p : b \leq q) \land  (\forall c<q  \exists  c<d<q : p \leq d).$$

Both functors recover the chronological relation on causally simple spacetimes.

\bigskip

Furthermore, there is a functorial definition of a topology $\tau_+$ for ordered spaces (going back to Beem and defined in generality in \cite{oM:Compl}), recovering the manifold topology for causally simple spacetimes and displaying several desirable properties, see \cite{oM:Compl}. 

For later use, we denote by $\CPOM_c$ the subcategory with those objects $(X, \leq, \mu)$ for which $ \mu (J(p,q)) \neq 0 $ for all $p,q \in X$ with $p \beta(\leq) q $ and such that $ (p,q) \mapsto \mu (J(p,q))$ is continuous w.r.t. $\tau_+$.
   
There are at least three possibilities to induce extended  (i.e., admitting the value $\infty$) pseudometrics on reasonably large subsets of ${\rm ALL^I}$ or $ {\rm POM^I}$.

 The first two approaches use functors from ${\rm ALL^I}$ resp. $ {\rm POM^I}$ to the category of metric spaces via a map $ \Phi_{f}$ (depending on $f \in {\rm LBM} (\R) $, being the set of locally bounded measurable real functions on $\R $) from an object $X$ to the space ${\rm AE}(X)$ of Borel-almost everywhere defined real functions on $X$, then defining ($p \in [1; \infty]$ fixed) a map $\Phi_{f,p}$ taking the class of $(X, \s) \in {\rm Obj} (\CALL)$ to the class of $ (X, d_{f,p})$ by
 
 \bea
 \label{PhiDef}
  \Phi_f : X \ni x \mapsto f \circ \s_x, \qquad  d_{f,p}(x,y) := \Phi_{f,p} (\s) := ((\Phi_f)^* d_{L^p}) (x,y)  = \vert f \ci \s_x - f \ci \s_y |_{L^p(X)}  \in [0; \infty ]
  \eea
 
 and $\s_x := \s (x, \cdot) $. On the metric side we have the Gromov-Hausdorff distance $d_{{\rm GH}}^+$ between isometry classes which we can pull back, inducing a map $\Phi_{f,p}^* d^+ _{{\rm GH}}: {\rm ALL^I} \times {\rm ALL^I} \rightarrow \R\cup \{ \infty \}$ resp. $\Phi_{f,p}^* d^+ _{{\rm GH}}: {\rm POM^I} \times {\rm POM^I} \rightarrow \R\cup \{ \infty \}$. In either of the two categories, we want to restrict ourselves to the quantities given there, without borrowing structures from the other category, because ultimately we want to show continuity of the (push-down of the) functor w.r.t. the respective Gromov-Hausdorff topology, which means that, on $\rm ALL^I$, we are only allowed to use the distance, limiting us to the choice\footnote{For brevity we denote by $|\cdot |_{\infty} $ the norm induced by the supremum and not by the essential supremum, whenever no canonical measure is part of the data.} $p = \infty$, whereas on $\rm POM^I$, we are only allowed to use the measure and the order, which implies $f := g \ci \sgn$ for some function $g$. Let us pursue the case of $\rm POM^I$ first. The calculations should be compared to the proof of Th. 1 (ii) in \cite{oM:LorFin}.

\bigskip

For $x \in X$, let $\s_x^+$ resp. $\s_x^-$ denote the positive resp. negative part of $\s_x$. For $r \in [-1;1] $ we define

$$ F_r:= - (\frac{1}{2} - \frac{r}{2})\chi_{(- \infty; 0)} + (\frac{1}{2} + \frac{r}{2})\chi_{(0;  \infty)} : \R \mapsto \R  , \ D_r := d_{F_r,2}, $$

then $D_r$ interpolates between the past metric (taking into account only the past cones) $D_{-1} $ with 

$$D_{-1} (x,y) := \vert \vert  \chi_{(-\infty; 0)}\ci \s_x -  \chi_{(-\infty; 0)} \ci  \s_y \vert \vert_{L^2 (X)} = \sqrt{\mu (J^-(x) \triangle J^-(y))} $$

 for $F_{-1} = (1- \theta_0)  $ and the future metric $D_1$ (taking into account only the future cones) for $F_1 = \theta_0 $, passing through $D_0 (x,y) = \frac{1}{2} \big || \sgn \s_x - \sgn \s_y \big||_{L^2}$. 
 Whereas $D_{\pm 1}^g$ is a metric on $ X \setminus \partial^\pm X$, it vanishes on $\partial^\pm X \times \partial^\pm X$ (which may be empty).  For fixed $p,q \in X$, let $V$ be the $2$-dimensional linear subspace ${\rm span} (u^+, u^-)$ spanned by the vectors $u^+ := \sgn \s_p^+  - \sgn \s_q^+$ and $u^- := \sgn \s_p^- - \sgn \s_q^-$. Then the $L^2$ scalar product on $V$ is uniquely given by the corresponding quadratic form on three vectors any two of which are non-collinear, thus given by $D_{-1/2} (p,q), D_0(p,q)$ and $D_{1/2} (p,q)$. In other words, the datum of those three recovers the whole family $D_r$. We are interested primarily in $D_{\pm 1} (p,q) = || u^\pm ||^2$. Elementary calculations reveal\footnote{The constants in Eq. \ref{RandRekuperation} are correct and differ from the incorrect ones appearing in \cite{oM:LorFin}, however, this correction does not affect any further statement or proof made in \cite{oM:LorFin}.}
 
 \bea
 \label{RandRekuperation}
   || u^\pm ||^2 =  || D_{\mp 1/2} (p,q) ||^2 + 3 || D_{\pm 1/2} (p,q) ||^2 - 3 || D_0(p,q) ||^2 . 
 \eea
 
 Now we want to recover the causal structure. For $p,q \in X$ we define $\kappa_p^\pm := \sgn (\s_p^\pm)$ and consider the following expressions only depending on the metrics $D_r$: 

$$ (D_{\pm 1 } (p,q) )^2     = \langle \kappa_p^\pm - \kappa_q^\pm , \kappa_p^\pm - \kappa_q^\pm \rangle_{L^2}   $$

$$ (D_0 (p,q))^2 = \frac{1}{4}\langle \kappa_p^+ - \kappa_p^- - \kappa_q^+ + \kappa_q^- , \kappa_p^+ - \kappa_p^- - \kappa_q^+ + \kappa_q^- \rangle_{L^2} , $$

we  calculate

\bea
\label{HarvestDiamonds}
   4 (D_0 (p,q))^2 -  (D_{-1 } (p,q) )^2 - (D_1 (p,q))^2 = 2 \big( \langle \kappa_q^+ , \kappa_p^- \rangle_{L^2} + \langle \kappa_p^- , \kappa_q^+ \rangle_{L^2}  - \underbrace{\langle  \kappa_p^+, \kappa_p^- \rangle_{L^2}}_{=0}  -  \underbrace{\langle  \kappa_q^+, \kappa_q^- \rangle_{L^2}}_{=0}   \big) ,   
   \eea   

where the last two terms vanish due to causality of $X$. This implies

\bea
 \label{DetectCausalRelation}
 (D_0 (p,q))^2 -  (D_{-1 } (p,q) )^2 - (D_1 (p,q))^2  \neq 0 \Leftrightarrow p \ll q \lor q \ll p .
 \eea

as $ p \ll q \Leftrightarrow \langle \chi_{(0;\infty)} \ci \s_p , \chi_{(-\infty; 0)} \s_q^- \rangle_{L^2} \neq 0$. The distinction between the two relevant cases on the right-hand side of Eq. \ref{DetectCausalRelation} can be made as follows: The above technique allows to determine the causal curves $c$ (without a direction). If $c: \R\rightarrow X$ is a future curve then $s \mapsto \lim_{t \rightarrow - \infty } D_1 (c(t), c(s))  = \mu (J^- (c(s)))$ is monotonically increasing in $s$, in the time-dual case it is monotonically decreasing. We define

$$\Phi^\times (X, \s) := (X, D_{-1/2}) \sqcup (X, D_0) \sqcup (X, D_{1/2})$$

 and $d_{{\rm GH}}^\times := (\Phi^\times)^* d_ {{\rm GH}}$ as the associated Gromov-Hausdorff metric on ${\rm POM^I}$, whose definition via distortion we now reconsider. 

A relation $R \subset X \times Y$ is called a {\bf correspondence} from $X$ to $Y$ if $\pr_1 R = X $ and $\pr_2 R = Y$, and the set of all correspondences between $X$ and $Y$ is denoted by $\Corr (X,Y)$.

\bea
\label{Corr1}
 d_{{\rm GH}}  ((X,d_X), (Y,d_Y)) := \frac{1}{2} \inf \{ \dist (R) \vert R \in \Corr(X,Y)  \} , 
 \eea

\V where for an arbitrary relation from $(X,d_X)$ to $(Y,d_Y)$ we define

\bea
\label{Corr2}
\dist (R) :=  \sup \{ \vert d_X(x,x') - d_Y(y,y') \vert : (x,y) , (x', y') \in R \}.
\eea

If we apply this to metric spaces $(X,d_X) $, $(Y,d_Y) $, then $d_{{\rm GH}}$ is indeed the Gromov-Hausdorff metric.

{\bf Remark.} The datum of $D_1$ alone is sufficient to reconstruct the (undirected) causal
structure, as it is easy to see that under the mild assumptions that $X$ is distinguishing and has open time cones and that any open set is of positive measure, a curve is $D_1$-geodesic if and only if it is causal. The author was not able to either construct an example of noninjectivity of $D_1$ on spacetimes or else proving injectivity. For our purposes, the construction above with $d_{{\rm GH}}^\times$ will prove sufficient. Moreover, $D_0$ could be used to topologize the elements of $\CPOM$. A reasonable condition would then be $\mu (I(p,q)) \neq 0$ for $p \ll q $ where $\ll$ is defined from $\leq $ via $\beta$ as below. Here, we will not pursue this approach further.  

\bigskip

Let ${\rm POM^I_{fv}}$ be the subset of ${\rm POM^I}$ of all those classes s.t. the measure of future and past cones is finite.

\begin{Theorem}
 $d^\times_{{\rm GH}} $ is an extended pseudometric on ${\rm POM^I}$ and a metric on ${\rm POM^I_{fv}}$.	Moreover, the volume of causal diamonds depends uniformly continuously on the $d^\times$-isometry class in the following sense: For each $\e >0$ there is a $\de >0$ such that for every correspondence $C: X \rightarrow Y $ between objects of $ \CPOM$ with $\dist^\times (C) < \de $ we have

$$ \forall (x_1, y_1), (x_2, y_2) \in C : | \mu_X (J_X(x_1,x_2)) - \mu_Y(J_Y (y_1, y_2))  | < \e . $$ 
\end{Theorem}

\V{\bf Proof.} The first statement is obvious, as the pull-back of a metric (here: the Gromov-Hausdorff metric) is a pseudometric. For the second statement, we only have to show injectivity of $\Phi_r^\times$ for $N= 3$, which has been shown in the preceding paragraph.
   The last statement follows directly from Eqs \ref{RandRekuperation} and \ref{HarvestDiamonds}. \hfill \qed

\bigskip

How to use the freedom of post-composition with some function $f$ in Eq. \ref{PhiDef} reasonably? In \cite{oM:LorFin} it has been shown that the metric defined by Noldus in \cite{Noldus} --- which in our terminology corresponds to $d_{g,| \cdot | , \infty}$ --- cannot be made a length metric by applying the "intrinsification" $\la$ to it\footnote{denoted by $\hat{\cdot}$ in \cite{BBI} and defined as usual by defining, for $c: I \rightarrow X$, $\ell (c) := \sup \{ \sum_{k=1}^{N-1} d(c(T_i, T_{i+1}) ) | T \in Y_N (I) \}$ where $Y_n$ is the set of ordered $N$-tuples in $I = [a;b]$ with $T_0 = a$ and $T_N = b$, and 
	$\la (d) (p,q) := \inf \{ \ell (c) | c: p \leadsto q {\rm \ continuous}\}$}. For $n \in \N$ we define $f_n(x) := \sgn (x) \cdot |x|^n$. For $n \geq 2$ and  $\check{d}_{n}:=  d_{g,f_n, \infty}$, the extended length metric $d_n := \la(\check{d}_{n}) \geq \check{d}_n$ is a length metric (i.e., finite) for $(X,g)$ a Cauchy slab (see Th. \ref{NoldusP}). Some properties are obvious: If $(X, \s)$ is {\bf distinguishing}, i.e., if for all $x ,y \in X$, $ x \neq y $ implies $I^\pm (x) \neq I^\pm (y) $, then $d_n$ is a metric on $X$ w.r.t. which $\s ^n$, and thus $\s$, is continuous, more precisely: $| \s^n (a,b) - \s^n (c,d) | \leq \check{d}_n (a,c) + \check{d}_n (b,d)$. Moreover, $|\s^n(x,y)| \leq \check{d}_n (x,y) $. However, the topology generated by $d_n$ differs in general from the topology $\tau_+$ defined in \cite{oM:Compl} (of course, if $\s$ is continuous w.r.t. $\tau_+$, as it is the case for $(X, \s)$ being globally hyperbolic, both topologies coincide). The functor $(X, \s) \mapsto (X, d_n) $ does not allow for a reconstruction of the Lorentzian data from the data on the metric side: Denoting by ${\rm ALL^I_{fd}}$ the subset of those isomorphism classes of compact almost Lorentzian pre-length spaces whose Lorentzian distance is finite, we obtain:

\begin{Theorem}
	\label{dplus}
	Let $f \in {\rm LBM}(\R)$. Then $({\rm ALL^I},\Phi_{f, \infty}^* d^+_{{\rm GH}})$ is an extended pseudometric space, and	
	
	$({\rm ALL^I_{fd}}, \Phi_{f, \infty}^* d^+_{{\rm GH}}|_{{\rm ALL^I_{fd}}\times {\rm ALL^I_{fd}}})$ is a pseudometric space but in general not a metric space, as in general, even in the case that $f: \R \rightarrow \R$ is a diffeomorphism, $\Phi_{f, \infty} $ is not injective. 
\end{Theorem}

\V{\bf Proof.} Let $ (X,\s)$ be defined by $X := [-1;1] \times \R$ (for $b \geq 1$) and $\s$ being the signed Lorentzian distance of the Lorentzian metric $k^2 \cdot (-dx_0^2 + dx_1^2)$, where $k \in C^{\infty} (X, (0;1])$ in such a way that there are no lines from $\partial^- X $ to $\partial^+X$ through $0$. For example, for $a \in (0; 1/4)$ we can choose $k_a (x)=1$ for all $x \in X \setminus (J(\{ 0 \} \times (-a;a)) \cap x_0^{-1} ((-1/4;1/4)) $ and $k_a (x)= a$ for all $x \in J((0,0)) \cap x_0^{-1} ((-1/4 +a; 1/4 -a )) $). Then, indeed, as each maximally extended line begins at $\partial^- X = \R\times \{ -1 \}  $ and ends at $\partial^+ X = \R\times \{ 1 \}  $, the image of no one can contain $(0,0)$: Any curve the image of which does has length $< 3 +5a$, but for each $p \in \partial^- X$ and each $q \in \partial^+ X$ there is a future curve from $p$ to $q$ of length greater or equal to the piecewise affine future curve $c$ from $(-1,0)$ to $(1,0)$ via $(-1/4,-1/4- a)$ and $(1/4, -1/4+a)$, which has length tending, with $a \rightarrow 0$, to $ 1/2 + \sqrt{8} >3$. That is, for sufficiently small $a$ there is no line through $(0,0)$. To construct another Lorentzian length $\tilde{\s}$ on $X$, we choose $f := \chi_{[0; \infty)}$. Now, for each two $p,q \in X$, the supremum of $\vert f \ci \s_p - f \ci \s_q \vert $ is attained on the boundary $\partial X$. We construct a Lorentzian length $\tilde{\s}$ on $X$ such that for each $(x,y) \in \leq$ we have a line from $x $ through $y$ to $\partial^+ X$, and such that for $A:= X \times \partial^+ X \cup \partial^- X \times X$ we have $\tilde{\s} \vert_A = \s \vert_A $, and thus $\Phi(\tilde \s) = \Phi (\s)$. It is easily seen via the inverse triangle inequality (which is sharp along lines) that the unique such distance $\tilde{\s}$ satisfies $\tilde{\s} (x,y) := \inf \{ \s (x,z) - \s (y,z) \vert z \in \partial^+ X \cap J^+(y) \} $ for all $x \leq y$ (for $y \leq x$ we extend by antisymmetry, on the complement by $0$). Let $A(X):= (X, \tilde{\s})$, then $ d_{{\rm GH}} (\Phi_f (X) , \Phi_f (A(X))) = 0$. Indeed, $A(X)$ is a pre-length space: For $x \ll y \ll z$,

\bean
 \tilde{\s} (x,y) + \tilde{\s} (y,z) &=& \inf \{ \s ( x,z) - \s (y,z ) \vert z \in J^+ (y) \cap \partial^+(y) \} + \inf \{ \s ( x,z) - \s (y,z ) \vert z \in J^+ (y) \cap \partial^+(y) \}\\
 &\leq& \inf \{ \s ( x,z) - \s (u,z ) \vert z \in J^+ (u) \cap \partial^+(y) \} = \tilde{\s} (x,u),
 \eean

and by the existence of the lines, $(X, \tilde{\s})$ is even a Lorentzian length space: For each two $x,y \in (X, \tilde{\s})$ there is a $\tilde\s$-maximal curve from $x$ to $\partial^+X$ through $y$. This property would be preserved by any isomorphism of $\CALL$, in particular for $y=0$ and $x \in \partial^- X$, thus $X$ and $A(X)$ are not isomorphic. \hfill \qed

\bigskip

As a third approach, we can mimic the definition of $d_{{\rm GH}}$ via correspondences and apply them directly to Lorentzian distance functions instead of metrics. Indeed, the definition of Gromov-Hausdorff distance via distortion by Eqs \ref{Corr1} and \ref{Corr2} equally well applies to $d_X,d_Y$ being signed Lorentzian distance functions, where we denote the assignment by $\dist^- $ and $d_{{\rm GH}}^-$ for better distinction. Let 

$${\rm tdiam} (X) := \inf \{ \sup \{ \s (x,y) \vert x \in \partial^- X \} \vert y \in \partial^+ X \} = \inf \{ \sup \{ \s (x,y) \vert y \in \partial^+ X \} \vert x \in \partial^- X \}$$

for $X \in {\rm Obj} (\CALL)$ and 

$$ {\rm tdiam} (X) := \inf \{ \sup \{ \mu ( I(x,y)) \vert x \in \partial^- X \} \vert y \in \partial^+ X \} $$  

for $X \in {\rm Obj} (\CPOM)$. 

Let ${\rm ALL}^I_{c}  \subset {\rm ALL^I_{fd}}$ resp. ${\rm CPOM}^I_c$ the subset of isomorphism classes of past and future compact g.h. objects $X \in {\rm Obj} (\CALL)$ resp. $X \in {\rm Obj} (\CPOM)$) with $\tdiam (X) >c$.

We equip the set ${\rm ALL^I}$ (resp. ${\rm ALL^I_{c}}  \subset {\rm ALL^I_{fd}}$) of isomorphism classes of objects (resp. past and future compact g.h. objects $X$ with $\tdiam (X) >c$) of $\CALL$ with the assignment $d^-_{{\rm GH}}$ and obtain:

\begin{Theorem}
	\label{ALLGH}
	$({\rm ALL^I},d^-_{{\rm GH}})$ is an extended pseudometric space, and	$({\rm ALL^I_{c}},d^-_{{\rm GH}})$ is a metric space. Moreover, $\Phi_{{\rm Id}}^* d^+_{{\rm GH}} \leq 2  d^-_{{\rm GH}}$. 
\end{Theorem}

\V{\bf Proof.} Nonnegativity and symmetry are obvious, and each isomorphism $ I: (X,d_X) \mapsto (Y,d_Y)$ is a correlation with $\dist^-(I) =0$. For the triangle inequality consider for two correspondences $ R_1 \in \Corr(X,Y)$ and $R_2(Y,Z) \in \Corr (Y,Z)$ their composition $R_2 \ci R_1 \in \Corr(X,Z)$. Then

\begin{align*}
	\dist^- (R_2 \ci R_1) &= \sup \{  \vert d_X(x,x') - d_Z(z,z') \vert : (x,z) , (x', z') \in R_2 \ci R_1 \}\\
	&= \sup \{  \vert d_X(x,x') - d_Y (y,y') + d_Y(y,y') - d_Z(z,z') \vert : (x,y) , (x', y') \in R_1 \land (y, z) , ( y', z') \in R_2 \}\\
	&\leq \sup \{  \vert d_X(x,x') - d_Y (y,y') \vert + \vert d_Y(y,y') - d_Z(z,z') \vert : (x,y) , (x', y') \in R_1 \land (y, z) , ( y', z') \in R_2 \}\\
	&\leq \sup \{  \vert d_X(x,x') - d_Y (y,y') \vert : (x,y) , (x', y' ) \in R_1 \} \\
	&   + \sup \{ \vert d_Y(\tilde{y},\tilde{y}') - d_Z (z,z') \vert : (\tilde{y},z) , (\tilde{y}', z' ) \in R_2 \} = \dist^- (R_1) + \dist^- (R_2)
\end{align*}

Now let $\e >0$ be given, let $d_{{\rm GH}}(X,Y) =: r$ and $d_{{\rm GH}}(Y,Z) =: s$, then there are correspondences $R_1: X \rightarrow Y$ and $R_2: Y \rightarrow Z$ with $ \dist^- (R_1) < r + \e/2$ and $\dist^- (R_2) < s + \e/2$ and thus $\dist^- (R_2 \ci R_1 ) < r +s + \e$, therefore $d_{{\rm GH}} (X,Z) \leq d_{{\rm GH}} (X,Y) + d_{{\rm GH}} (Y,Z) $. 
To see $d_{{\rm GH}}(X,Y) = 0 \Rightarrow X=Y$ and completeness, we first note $\dist^- (R) = \dist^- (\cl (R, X \times Y)) $. Obviously, every correspondence $f: X \rightarrow Y$ between distinguishing ordered measure spaces $X$ and $Y$ with $\dist^- (f) = 0 $ is a bijective map and an isomorphism of $\CALL$. Now assume $d^-_{{\rm GH}} (X,Y) =0$. Then there is a sequence $a$ in $\Corr(X,Y) $ with $\dist^-(a_n) \rightarrow_{n \rightarrow \infty} 0$. Let $\Sigma$ be a dense countable set in $X$, which exists due to compactness and metrizability (by $d^+$) of $X$. Compactness of $Y$ and Cantor's diagonal procedure ensures existence of a subsequence $ b$ of $a: n \mapsto a_n$ such that for every $s \in \Sigma$ some points in $b_n(s)$ converge to some $\tilde{b}_\infty (s) \in Y$, so $\tilde{b}_\infty : \Sigma \rightarrow Y$. Let $b_\infty := \cl (\tilde{b}_\infty, X \times Y)$, then $b_\infty$ is left-total, as $\Sigma$ is dense in $X$, and it is right-total, as by $\dist^+ < 2 \dist^- $, the sets $b_n(\Sigma)$ are dense in $Y$. Thus $b_\infty$ is a correlation. For all $(p,r), (q,s) \in b_n$, we get $ \vert \s_Y (r, s) - \s_X (p,q)  \vert \leq \dist^- (b_n) \rightarrow_{n \rightarrow \infty} 0$, thus $\s (b_\infty (p), b_\infty (q) ) = \s (p,q) $ and $b_\infty$ is an $\CALL$-isomorphism. 

Now we show the estimate between the metrics: For every correspondence $R: X \rightarrow Y$ we have 

\bean
\dist^+(R) &=& \sup_{(x_1, y_1) , (x_2,y_2) \in R } \big( \big\vert \sup_{x_3 \in X} \vert \s (x_1, x_3) - \s (x_2, x_3) \vert - \sup_{y_3 \in Y} \vert \s (y_1, y_3) - \s (y_2, y_3) \big\vert \big)
\eean

Supposing w.l.o.g. that $ \sup_{x_3 \in X} \vert \s (x_1, x_3) - \s (x_2, x_3) \vert \geq \sup_{y_3 \in Y} \vert \s (y_1, y_3) - \s (y_2, y_3) |$, we get

\bean
\dist^+(R) &\leq& \sup_{(x_1, y_1) , (x_2,y_2) \in R } \big(  \sup_{x_3 \in X} \vert \s (x_1, x_3) - \s (x_2, x_3) \vert - \sup_{y_3 \in Y} \vert \s (y_1, y_3) - \s (y_2, y_3)  \big)\\
&\leq& \sup_{(x_1, y_1) , (x_2,y_2) \in R } \big(  \sup_{x_3 \in X, y \in R(x_3)} ( \vert \s (x_1, x_3) - \s (x_2, x_3) \vert -  \vert \s (y_1, y_3) - \s (y_2, y_3) ) \big)\\
&\leq&  \sup_{(x_1, y_1) , (x_2,y_2) , (x_3,y_3 )\in R } \big( | \s (x_1, x_3)  - \s (y_1, y_3)| + | \s (x_2, x_3)  - \s (y_2, y_3)|  \big)\\
&=& 2 \dist^- (R)
\eean

(where the inequality in the second last line is an application of the general rule $\big| | A-B | - | C - D| \big| \leq \big| (A-B) - (C-D) \big| \leq | A-C| + |B-D|  $ for all real numbers $A,B,C,D$). Thus $d_{{\rm GH}}^+ \leq 2 d_{{\rm GH}}^-$. \hfill \qed

\bigskip

As an example, we let $r \in (0; \infty]$ and consider $ L_{n,r} := \frac{1}{n} \Z^{1,1} \cap x_0^{-1} (-r;r) \subset \R^{1,1} \cap x_0^{-1} (-r;r)=: L_{\infty,r}$. Then considering $R_n \in \Corr(L_{n,r}, L_{\infty,r}) $ with $xR_ny : \Leftrightarrow y \in B^{\R^2} (x, 2/n)$, we see $L_n \rightarrow_{n \rightarrow \infty} L_\infty$ for $r < \infty$ but not for $r = \infty$ in the (non-Hausdorff) topology generated by the pseudometric $\Phi^*_{f, \infty}d_{{\rm GH}}$ for every continuous bijective function $f$. For further example, consider the second last paragraph of this section.

\bigskip

McCann and S\"amann \cite{rMcS} construct a natural measure on Lorentzian pre-length spaces (and so a functor $\MCS: \CALL \rightarrow \CPOM $), essentially equivalent (see \cite{oM:LorFunct}) to the following: For $A \subset X$, $N>0$ we define (keeping in mind that $I(x,y) $ is open with compact causally convex closure for all $x,y \in X$ if $X$ is g.h.):

$CC_{\delta} (A) := \{ (p,q) \in (X^{\N})^2 \vert A \subset \bigcup_{k=1}^{\infty} J(p(k), q(k)) \land \diam (I (p(k), q(k))) < \delta \forall k \in \N  \},$

(where $\diam$ is w.r.t. the metric $d_{U_k}$ for $U_k:= I(p_k,q_k)$),

$$ \la_N (D) := \omega (N) \cdot \sum_{k=1}^\infty \s (p(k),q(k)))^N \forall D \in CC_{\delta} (A), \ {\rm where} \  \omega (N):= \frac{\pi^{\frac{N-1}{2}}}{N \cdot \Gamma (\frac{N+1}{2}) \cdot 2^{N-1}}, \ {\rm and} $$

$$ \mu_{N, \delta} (A) := \inf \{ \la_N(D) \vert D \in CC_{\delta} (A)\}, \qquad  \mu_N (A) := \lim_{\delta \rightarrow 0} \mu_N (\delta) = \sup_{\delta >0} \mu_N(\delta),  $$

\V The outer measure $\mu_N$ induces a unique measure $\MCS (\s)$ on the Borel subsets. 

Conversely, in \cite{oM:LorFunct}, the author defined, for an ordered measure space $(X, \leq, \mu) \in {\rm Obj} ( \CPOM_c )$ and $b \in X$, the {\bf dimension} $\Dm (b)$ {\bf of $X$ at $b$} (which recovers the dimension $n $ at every $b \in X$ if $X$ is a spacetime) by

$$ \Dm(b) := \limsup_{I^+ (b) \ni c \rightarrow b \leftarrow c \in I^-(b)} (-\log_2 (\Phi (a,c)) +1) \in [1; \infty].$$

where, for $a,c \in X$ with $a \ll c$ (and thus $\mu (J(a,c)) \neq 0$ by definition of $\CPOM_c$)

$$ \Phi(a,c)  :=  \sup \Big\{ \frac{\mu(J(a,b) ) + \mu(J(b,c) )}{\mu (J(a,c)) } \Big\vert b \in J(a,c)  \land \mu (J(a,b)) = \mu (J(b,c))  \Big\} \in [1; \infty] ,$$

and defined the Lorentzian length $\ell(c) \in [0; \infty]$ of a timelike curve $c: I \rightarrow X$ by

$$ \ell(c) = \inf \Big\{ \sum_{k=1}^N \sqrt[\Dm(p_k)]{\frac{\Dm(p_k) \Gamma (\frac{\Dm(p_k)-1}{2} +1)}{\pi^{\frac{\Dm(p_k)-1}{2}}} \cdot \mu (J(p_k, p_{k+1}))}  \Big\vert \{t_0 ,... t_N \} {\rm partition \ of \ } I {\rm \ with \ }  p(n) := c(t_n) \Big\}, $$

and finally put $\mathcal{F} (X, \leq, \mu) := (X, K^-(\ell))$ where $K^- (\ell ) (x,y) := \sup \{ \ell (c) \vert c: x \leadsto y {\rm \ causal}\}$.

\bigskip
 
We easily see 
 
\begin{Theorem}
Let $c>0$. Both $ {\rm POM}^I_c$ and ${\rm ALL}^I_c$ are closed, and	$\mathcal{F} (\CPOM_c) \subset \CALL_c$. \qed
\end{Theorem}

\bigskip

Let us give a statement of equivalence of the topologies generated by the two Gromov-Hausdorff metrics.

Let ${\rm POM^I}_{n,c}$ resp. ${\rm POM^I}_{\leq n, c}$ be the space of compact elements of ${\rm POM^I}_{c}$ of constant dimension $n$ resp. of local dimension $\leq n$ and let ${\rm ALL^I}_{n,c}$ resp. ${\rm ALL}^I_{\leq n, c} $ be the space of compact elements of ${\rm ALL^I}_c$ of constant dimension $n$ resp. of local dimension $\leq n$. 
Moreover, let, for $\e \in (0; {\rm tdiam} (X)/2)$ , we define $K_\e:  \CALL \rightarrow \CALL$ preserving (future/past) compactness, connectedness and global hyperbolicity by

$$K_\e (X) := \{ p \in X \vert \sup \{  \s (p,y) \vert y \in \partial^+ X  \} \geq \e  \land \sup \{ \s (x,p) \vert x \in \partial^- X  \} \geq \e   \} .$$

\begin{Theorem}
	\label{ContinuousFunctor}	
Let $\e >0$ and $n \in \N$. Then $ \mathcal{MCS} \ci K_\e: ({\rm ALL^I}_n, d_{{\rm GH}}^-) \rightarrow ({\rm POM^I}, d_{{\rm GH}}^{\times})$ and 

$ K_\e \ci \mathcal{F} : ({\rm GOM^I}_n, d_{{\rm GH}}^\times) \mapsto ({\rm ALL^I}, d_{{\rm GH}}^-)$ are continuous.
\end{Theorem}

\V{\bf Proof.} For the first assertion, let $X_\infty \in {\rm ALL^I}_n$ be given, we want to show continuity of $\mathcal{MCS} \ci K_\e $ at $X_\infty$. Let $n \mapsto (C_n: X_n \rightarrow X_\infty)$ be a sequence of correspondences between objects of $\CALL$ with $ \dist^- (C_n ) \rightarrow_{n \rightarrow \infty} 0$. It is sufficient to prove that $\dist^\times (\mathcal{MCS} (K_\e(C_n))) \rightarrow_{n \rightarrow \infty} 0$. We want to show that the distortion w.r.t. the outer measure $\mu_N$ converges to $0$. Assume the opposite, then there is $\e >0 $ and points $p_n, q_n \in X_n$ as well as $ p_\infty , q_\infty \in X_\infty$ with $(p_n, p_\infty ) , (q_n,q_\infty) \in C_n$ and $ \vert d^\times (p_n, q_n ) - d^\times (p_\infty, q_\infty) \vert > \e $. On the other hand, we have $  d^\times (p_n, q_n ) = \vert \s_{p_n} - \s_{q_n} \vert_{L^2(X_n)}$ and $  d^\times (p_\infty, q_\infty ) = \vert \s_{p_\infty} - \s_{q_\infty} \vert_{L^2(X_\infty)}$. As we know from the proof of Th. \ref{ALLGH} that $\dist^+ (C_n) \leq 2 \dist^- (C_n)$, the respective integrands (in the definition of the $L^2$ norm) become arbitrarily close in the supremum norm, thus the only thing that remains to be shown to complete the proof by contradiction is that the measures of the $X_n$ are uniformly bounded. To this aim, let $\de >0$ be fixed. Let, for $\e >0$, $I_\e (x,y) := \{p \in X \vert \s (x, p ), \s (p,y ) > \e\}$,

$$CC_{\delta}^\e (A) := \{ (p,q) \in (X^{\N})^2 \vert A \subset \bigcup_{k=1}^{\infty} I_\e(p(k), q(k)) \land \diam (I_\e (p(k), q(k))) < \delta \forall k \in \N  \}.$$

Let $A:= K_\e (X_\infty)$ and $D_\infty \in CC_{\de}^{3 \e} (A)$ with $\la_n (D_\infty) > \mu_{n, \de} (A) - \theta$. By compactness of $A$, there is a finite subcover $D_\infty^0$. Then we choose, for every $(p,q) \in D^0_\infty$, an element $(p^{(m)},q^{(m)})$ with $(p^{(m)},p), (q^{(m)},q) \in C_m$ and obtain a finite cover $ D_m  \in  CC_{\de + 2 \e} (X_m)  $ of cardinality $N$. Indeed, $\dist^+ (C_m) < 2 \e$ ensures that each $x \in X_n$ is contained in some $ I(p_k,q_k)$. If $ \s(p,q) \geq \e$, then $\s(p^{(m)}, q^{m()}) < 2 \s (p,q)$, whereas if $\s(p,q) \leq \e$, then $\s(p^{(m)}, q^{(m)}) \leq 3 \e$. Thus $\la_n (D_m) <  2^n \la_n (D_\infty) +  N (3\e)^n \leq 2^n (\mu_{n,\de} (A) - \theta) + N (3 \e)^n$.

For the second assertion, let $n \mapsto (C_n: Y_n \rightarrow Y_\infty)$ be a sequence of correspondences between objects of $\CPOM$ with $ \dist^\times (C_n ) \rightarrow_{n \rightarrow \infty} 0$. We want to show $\dist^- (\mathcal{F} (C_n)) \rightarrow_{n \rightarrow \infty} 0$. Assuming the opposite, there is $\e>0$ and a sequence of elements $x_n, y_n  \in X_n$ and $a_n , b_n \in X_\infty $ with $ (x_n, a_n), (y_n, b_n ) \in R_n$ and $| \s (x_n, y_n ) - \s (a_n, b_n) | > \e $. By compactness of $X_\infty$, the sequences $a$ and $b$ accumulate at  some $x_\infty, y_\infty \in X_\infty$. Let $r,s \in J(x_\infty, y_\infty)$ with $\s (r,s) > \s (x_\infty, y_\infty) - \e/3$. For large $n$, $a_n \in I^- (r)$ and $ b_n \in I^+(s)$. Let $c$ be a causal curve from $r$ to $s$ with $ \ell (c) > \s (x_\infty, y_\infty ) - \e/2  $, thus there is a finite causal chain $p$ from $a_n$ to $ b_n$ with $\ell(p) > \s (x_\infty, y_\infty ) - \e/2 $. Let $M := \# p$. We want to pull back the causal chain to a causal chain in $X_n$. Let $\de:= \min \{ d^\times (p_i, X \setminus J^-(p_{i+1})) \vert i \in \N_{M-2}\}$. Let $p_i^{(n)} \in C_n^{-1} (p_i)$. There is $n_0 \in \N$ with $\dist^\times (C_n) < \{ \frac{\e}{2M}, \de \}$ for all $n \geq n_0$. Then $p^{(n)}$ is a causal chain from $x_n$ to $y_n$ for all $n \geq n_0$ and $\s (x_\infty, y_\infty ) - 2 \e < \ell (p^{(n)}) \leq \s (x_n,y_n) $. This shows lower semicontinuity $\s(x_\infty, y_\infty) \leq \liminf_{n \rightarrow \infty} \s (x_n,y_n)$. For upper semicontinuity, choose a subsequence $n \mapsto \gamma_n$ of $n \mapsto c_n$ with $ \gamma_n \rightarrow_{n \rightarrow \infty}^{d^\times} \gamma_\infty $ and $ \tilde{\gamma}_n $ by dyadic convergence, i.e. successively choosing midpoints (there is a Lipschitz parametrization of each curve). Then $\tilde{\gamma}_n \rightarrow_{n \rightarrow \infty} \gamma_{\infty}$, which implies the claim by upper semicontinuity of $\ell$ (\cite{mKcS}, Prop. 3.17). \hfill \qed

\bigskip

{\bf Remark.} As the two metrics scale differently for $n \neq 1$, they are not Lipschitz equivalent.

\begin{Theorem}
	\label{InversesOnClosure}
The functors $\mathcal{F}$ and $\mathcal{MCS}$ are inverse to each other on the closure of the set of classes of g.h. spacetimes in ${\rm ALL^I}_n$ resp. ${\rm POM^I}_n$.	
\end{Theorem}

\V{\bf Proof.} Assume they are not inverse to each other on an element $X$ of the closure, then there is some $\e>0$ such that $ \mathcal{F} \ci \mathcal {MCS } \ci K_\e (X) \neq  K_\e  (X)$. But the last equality holds for spacetimes, the maps are continuous by Th. \ref{ContinuousFunctor}, so the equality holds for the closure, contradiction. \hfill \qed

\bigskip
 
Interestingly, despite the fact that $d_{{\rm GH}}^-$ is not a metric between metric spaces, there is an adapted Gromov precompactness result in this context. For a full and g.h. almost Lorentzian length space $X$ we could define an $\e$-net to be a subset $N \subset X$ such that for each $x \in X$ there are $p,q,y \in N$ with $x,y \in J(p,q)$ and $ \vert \vert \s_x - \s_y \vert \vert_{L^\infty (J(p,q))} < \e $, and then mimick BBI 7.4.9, 7.4.12, 7.4.15 and 7.5.1 and show completeness. However, in the following we will use another way. 
 Let ${\bf MET}$ the category of metric spaces. Then for each $\e >0$ and $N \in \N$ and each $A \subset {\rm MET^I} $ such that $\diam (X) \leq T$ for all $X \in A$ we get

\bigskip

{\em For every $X \in A$, there is an $\e$-net of cardinality $N$ in $X$ 
	
	$\Leftrightarrow $  there is a metric $d$ on $\N_N^*$ such that $d_{{\rm GH}} (X, (\N_N^*,d)) < \e/2$
	
	$\Leftrightarrow A \subset B_{d_{{\rm GH}}} ( A_{N,T+ \e}, \e/2)$ for the (compact) set $A(N,T+ \e) $ of isometry classes of metric spaces of cardinality $\leq N$ and diameter $\leq T + \e$.}  

\bigskip

Thus one could define: $X$ uniformly totally bounded $ : \Leftrightarrow \forall \e >0 \exists N \in \N \exists T>0 : A \subset B_{d_{{\rm GH}}} (A_{N,T}, \e)$. This is the definition most useful for our aims as it persists for non-metric Gromov-Hausdorff spaces. We first note an auxiliary theorem of independent interest:

\begin{Theorem}
	\label{PointwiseLimit}
Each pointwise limit of Lorentzian distance functions is a Lorentzian distance function.	
\end{Theorem}
 
 \V{\bf Proof.} This is obvious from the fact that the two conditions (antisymmetry and conditional inverse triangle inequality) each concern the $2$- or $3$-tuples of points in the respective set, their hypotheses are open conditions and their conclusions are closed conditions. \hfill \qed

\begin{Theorem}[Gromov precompactness for $d_{{\rm GH}}^-$]
Each $A \subset {\rm ALL^I_{c}}$ is precompact w.r.t. $d_{{\rm GH}}^-$ iff

\begin{enumerate}
	\item There is $T >0$ with $\diam^- (X) := \sup \{ |\s (x,y)| : x, y \in X \} < T \ \forall X \in A$;
	\item For all $\e >0$ there is $N (\e) \in \N$ such that for all $X \in A$ there is $ Y_\e \in {\rm ALL^I}$ with $ \# Y_\e \leq N(\e) $ and $d_{{\rm GH}}^- (X,Y_\e) < \e$. 
\end{enumerate}

\end{Theorem}

\V{\bf Remark.} The conditions imply compactness of each $X \in A$, as the objects of $\CALL_c$ are complete.
\V{\bf Proof.} We adapt the proof of \cite{BBI} for metric spaces to the present situation. Let $X: n \mapsto X_n$ be a sequence in $A$. Let $N(\e)$ as above always be chosen minimal.

Let $X \in A$ and let $Y_{1/m}$ satisfy the hypothesis (Item 2 of the theorem) with $\e = 1/m$. Let $E^{X,m} \in X^{N(1/k)}$ be in correspondence of distortion $\e$ with $ Y_{1/m}$. We define recursively a sequence $M: \N\rightarrow \N$ by $M(0) := 0$ and $M(k+1) := M(k) + N(\frac{1}{k+1}) $ for all $k \in \N$, and a sequence $a_X: \N \rightarrow X$ by 'joining $1/k$-nets': For $\mu (n) := \sup \{ k \in \N \vert M(k)  \leq n \}$, let $a_X (l) = E^{X, \mu (k)}_{l - \mu(k) +1}$.

Now we define a sequence of functions $\s_n: \N\times \N\rightarrow \R  $ by $ \s_n (u,v ) := \s_{X(n) } (a_{X(n)} (u) , a_{X(n)} (v))$. Those functions on $\N\times \N$ still satisfy both properties of a signed Lorentzian distance function. As $\s_n (u,v ) \leq T$ for all $n,u,v$, by Cantor's diagonal procedure we can find a subsequence $c= b \ci j $ of $b$ such that $\s_{j(n)} (u,v)  $ converges in $\R$ for all $u,v \in \N$. Let $\s_{\infty} : \N\times \N\rightarrow \R $ defined as a pointwise limit by $\s_{\infty} (u,v)  := \lim_{n \rightarrow \infty} \s_{j(n)} (u,v) $ for all $u,v \in \N$. Then by Theorem \ref{PointwiseLimit}, $\s_\infty$ is still a signed Lorentzian distance function. Now, $d^+ : \N \times \N \rightarrow [0, 2T] $ defined by $d^+ (m,n) := \sup \{ \vert  \s_\infty (m,k) - \s_{\infty} (n,k) \vert : k \in \N  \}$ (in accordance with the paragraph preceding Theorem \ref{dplus}) is a pseudo-metric on $\N$, and we define the canonical quotient map $q$ by the vanishing subset to get a metric space $(\tilde{X}_\infty, \tilde{d}^+)$. Still, $\s_\infty$ can be pushed down in a  well-defined way to a signed distance function $\tilde{\s}$ on the quotient by the very definition of $d^+$, and can, for the same reason, be canonically extended to a signed Lorentzian distance function $ \hat{\s}$ on the Cauchy $d^+$-completion $(X_\infty, D)$ of $(\tilde{X}_\infty, \tilde{d}^+)$ by assigning to two $x_\infty, y_\infty \in X_\infty) $ the distance $\hat{\s} (x_\infty,y_\infty) := \lim_{n \rightarrow \infty} \tilde{\s} (x(n),y(n))$ for two arbitrary sequences $x,y$ in $\tilde{X}_\infty$ converging to $x_\infty$ resp. $y_\infty$. Moreover, $D= d^+ (\hat{\s}) $. Now our aim is to show compactness of $(X_\infty,D)$ and $d_{{\rm GH}}^-$-convergence of $X$ to $(X_\infty,  \hat{\s})$.   
 
 As for compactness, let $S^{(k)}$ be the image of $q(\N_{N(k)})$ in the completion, where it is an $1/k$-net in the sense that the inclusion of $S^{(k)}$ can be extended to a correspondence of $\dist^- < \frac{1}{k}$. To see that this is true, recall that for $ S_n := a_{c(n)} (\N) $, the inclusion $I$ of $ a_{c(n)} (\N_{N(k)}) $ in $S_n$ has $d_{{\rm GH}}^- $-distance $\leq 1/n$ from $c(n)$. As $d_{{\rm GH}}^+ \leq 2 d_{{\rm GH}}^-$ by Theorem \ref{ALLGH}, we have $d_{{\rm GH}}^+ (c(n), A(N_k,T) )$. But $N_k $ does not depend on $n$. Thus for each $i \in \N$ there is $ j \in \N_{N_k}$ with $d^+(x_{i,n}, x_{j, n})$ for all $n \in \N$. Thus $d^+ (q(x_i) m q(x_j)) \leq 1/k$. Consequently, $S^{(n)}$ is an $1/k$-net (w.r.t. $d^+$) in $\tilde{X}$ and also in its completion $X$. Therefore the complete space $X$ has a finite $1/k$-net for all $k \in \N$ and is therefore compact. A bit more conceptually: As $S^{(k)} \in A(N,T)$, we get $\forall k \in \N\exists N \in \N\exists T>0: X \in B^{d_{{\rm GH}}^+} (A(N,T), 1/k)$.

Convergence is readily seen by the fact that $S_n^{k} \rightarrow^{d_{{\rm GH}}^-}_{k \rightarrow \infty} S_n$ as finite sets (the Lorentzian distances converge), but then by the above estimate the sequence also converges w.r.t. $d_{{\rm GH}}^+$. As this is true for every $k \in \N$, the required convergence holds by the familiar nets criterion for Gromov-Hausdorff convergence (see \cite{BBI}, Prop. 7.4.12, e.g.). 
\hfill \qed 

\bigskip

\V{\bf Remark.} If we want to show global hyperbolicity of a Gromov-Hausdorff limit $X$ of globally hyperbolic Lorentzian length spaces then by the a priori restriction to compact Lorentzian length spaces we only have to show causal simplicity (closedness of $J^\pm(x)$ for all $x \in X$), which is automatic by the choice of the topology $\tau_+$ as in \cite{oM:Compl}. 

\bigskip

Recall that each compact metric spaces can be approximated in $d_{{\rm GH}}^+$ arbitrarily close by finite metric spaces --- which are not length metric spaces. The Hausdorff dimension of a Gromov-Hausdorff limit of length spaces $X_n$ under a uniform lower curvature bound is $\leq \limsup_{n \rightarrow \infty} \dim (X_n)$.  This motivates the

{\bf Question:} Howe large is the closure of the set of spacetimes in the set of globally hyperbolic Lorentzian length spaces satisfying a lower curvature bound?

By Theorem \ref{InversesOnClosure}, an answer would shed some light on the question to which extent the functors $\mathcal{F}$ and $\mathcal{MCS}$ are inverse to each other on ${\rm ALL^I}_{\leq n} $ resp. $ {\rm POM^I}_{\leq n}$.

\bigskip

As it becomes more and more transparent that ordered measure spaces and almost Lorentzian length spaces are intimately related such that, under mild hypotheses, both structures can be used at a time, one can try to set up a version of analysis on metric measure spaces as initiated by Ambrosio, Gigli and many others (for an overview see \cite{nG1}, \cite{nG2}), in this new context. The theory of Lorentzian optimal transport as in \cite{Stefan}, \cite{MS} naturally defines Ricci curvature bounds.

\bigskip

To define a pointed version for the noncompact case, there are two possible lines of attack:

\begin{enumerate}
	\item Begin with Gromov's pointed metric as laid down in Herron \cite{dH} and try to reformulate it via correspondences.
	\item Use the Busemann metric instead of the Hausdorff metric and try to reformulate the corresponding Gromov-Busemann distance in terms of correspondences.
\end{enumerate}

We would need to replace balls being compact subsets exhausting the spaces (as the pointed Hausdorff distance is a true metric only on the set of isometry classes of proper metric spaces). A small consideration shows that there is no $SO(1,n)$-covariant map from $\R^{1,n}$ to the set of open precompact subsets of $\R^{1,n}$. So to define a natural metric we need at least a double puncture on each space, corresponding to the well-known {\em observer moduli space}, a quotient by those diffeomorphisms fixing a given tangent vector. 

\bigskip

For a Lorentzian length space $(X, \s)$ and a complete length space $(Y,d)$ we define their {\bf product} $(X \times Y, p(\s, d) )$ by $p(\s, d) ((x_1, y_1), (x_2, y_2)) := \pm \sqrt{\s^2 (x_1, x_2) - d^2 (y_1,y_2)} $ for $\pm (\s(x_1, x_2 ) > d(y_1,y_2)$ and $0$ otherwise. It is easy to see that $(X \times Y, p(\s,d))$ is a Lorentzian length space, and that we have $d_{GH}^- (X \times Y, A \times B) = \sqrt{d_{GH}^- (X,Y)^2 + d_{GH}^+ (A,B)^2}$ (as the natural product of correspondences satisfies the corresponding equation). Furthermore, we obtain $d_{GH}^- ((X, r \cdot \s), (X, \s)) = r \cdot {\rm tdiam} (X, \s)$ for all $r >0$ if $ {\rm tdiam} (X, \s)  <\infty$ in exactly the same manner as in the metric case, i.e. by calculating the distortion of the identity. This implies e.g. that the curve $r \mapsto ((\R, -dt^2) \times \mathbb{S}^n (\frac{r}{2 \pi}))  $ is geodesic, where $\mathbb{S}^n (r) = r \cdot \mathbb{S}^n (1)$ is the standard sphere of radius $r$ in Euclidean $\R^n$.

\bigskip

Only two days after the prepublication of the first version of this article on arXiv.org (whose Sec.2 differed from Sec.2 of the present article only by minor details and this last paragraph), Ettore Minguzzi and Stefan Suhr prepublished on arXiv.org the results of their independent work \cite{eMsS} on the same subject in a very recommendable article which contains other very interesting aspects, like stability under $d_{GH}^-$-convergence in the set of Lorentzian pre-length spaces (a) of the property of being a so-called bounded Lorentzian length space with $i_0$ (\cite{eMsS}, Th. 5.18) and (b) of timelike curvature bounds (\cite{eMsS}, Th. 6.7). 

\newpage

\section{Dimensions of ordered sets}
\label{Dim}

Notions of dimension are ubiquitous in the different branches of mathematics. All of them relate to cross products and reflect the existence of left or right invertible morphisms between objects. The author suggests the following general definition of dimension: Let $C$ be a category, a {\bf dimension} in $C$ is a map from ${\rm Obj} (C)$ into a totally ordered set\footnote{Often the totally ordered set is one of cardinalities. One can remain in the domain of set theory by just bounding the cardinality of the objects. Of course for this we need that the cardinality of the product of two infinite nonempty sets $A,B$ is equal to $\max \{ \# A, \# B \}$.} $D$ such that for all $A,B \in {\rm Obj} (C)$ we have

\bea
\label{Dim-Def}
\dim (A) < \dim(B) \Rightarrow \not \exists {\rm \ injective \ morphism \ } B \rightarrow A .
\eea

We call a dimension $\dim$ {\bf strong} iff it satisfies $\dim (A) < \dim(B) \Rightarrow \not \exists {\rm \ surjective \ morphism \ } A \rightarrow B$. 
If $C$ has a forgetful functor to the monoidal category of sets and maps with the Kuratowski cross product and $\dim (A \times B) = \dim (A) + \dim (B)$, we call $\dim$ {\bf monoidal}.

As it is to be expected by the form of the definition, dimensions are intimately related to Cantor-Bernstein theorems in the respective categories. Of course the prime example is $C_0 $ being the category of sets with the usual Kuratowski product and $\dim = \#$. Other examples (many of which are strong and monoidal):

\begin{enumerate}
	\item $C^K_1$ the category of $K$-vector spaces and $\dim_1^K (X)$ the cardinality of a Hamel basis of $X$,
	\item $C_2^K$ the category of $K$-Hilbert spaces and $\dim_2^K$ the cardinality of a Hilbert basis of $X$, 
	\item $C_3$ the category of topological spaces and $\dim_3$ being the small inductive, large inductive, or the covering dimension (all coinciding on the separable metrizable spaces),
	\item $C_4$ the category of metric spaces and Lipschitz maps with the Hausdorff dimension,
	\item $C_5$ the category of metric spaces and coarse maps and Gromov's asymptotic dimension,
	\item $C_6$ the category of manifolds and embeddings with the usual manifold dimension,
	\item $C_7$ the category of ordered sets and increasing maps with the product order on the cross product with several dimensions, e.g. with the Dushnik-Miller dimension.	 
\end{enumerate} 

In the last example, the {\bf product order} on an arbitrary Cartesian product $\Pi_{i \in I} C_i$ defined by $ x \leq y : \Leftrightarrow x_i \leq y_i \ \forall i \in I$. A map $f$ between ordered spaces $X$ and $Y$ is called {\bf order embedding} iff $ a \leq  b \Leftrightarrow f(a) \leq f(b)$. Most of the above notions of dimension have the virtue of being compatible, i.e., of coinciding on the intersections of the categories, the only two exceptions being firstly the well-known inequality between the Hamel dimension $\dim^K_1$ and the Hilbert dimension $\dim^K_2$ in infinite dimensions and secondly the inequality beween the Dushnik-Miller dimension on one hand and all other dimensions on the other hand, as we will see below. All dimensions above refine cardinality in the sense that $ \#A < \# B $ implies $\dim (A) < \dim (B)$ for each two objects $ A,B$.

Let $(X, \leq)$ be an ordered set.  The {\bf Dushnik-Miller dimension $\dim_{DM} (X)$ of $X$} is defined as

$$ \dim_{DM} (X) := \inf \{ \# \{ A \subset P( X \times X ) \} \vert \forall a \in A: a {\rm \ total \ order \ on \ } X {\rm  \ and \ } \leq = \bigcap A  \} .$$ 

Dushnik and Miller \cite{bDeM} showed that $\dim_{DM} (X)$ is well-defined and $> - \infty$ as there are subsets $A$ satisfying the conditions, and Ore \cite{oO} showed that $\dim_{DM} (X)$ is the smallest cardinal $I$ such that there is an order embedding $f$ from $X$ into the Cartesian product $ (\Pi_{i \in I} C_i, \leq_I)$ where the $ C_i$ are totally ordered sets and $\leq_I$ is the product order\footnote{for the nontrivial direction use the well-ordering theorem to choose a well-ordering $\leq$ on $I$ and consider, for each $i \in I$, the bijection $B_i$ commuting $ i $ and the minimal element $0$ leaving anything else unchanged, then let $\leq_j$ be the lexicographical order on $\Pi_{i \in I } C_i $ w.r.t. the total order $B_j^* \leq$ on $I$, which is again total, and the product order is the intersection $ \bigcap_{j \in I} \leq_j$}. It can easily be seen that $\R^n_{prod}$, which denotes $\R^n$ with the product order, has Dushnik-Miller dimension $n$, in particular $\dim_{DM}(\R^{1,1}) =2$. Moreover, the power set $P(S)$ of a set $S$ equipped with the order $ \leq = \subset$ has Dushnik-Miller dimension $\# S$. 

The authors of \cite{HBG} defined two other notions of dimension for ordered spaces. To give an account of those, let us first define some other notions. A {\bf utility on $X$} is a map $f: X \rightarrow \R^I$ for a set $I$ such that $\forall x ,y \in X: x \leq y \Leftrightarrow f_i(x) \leq f_i(y) \ \forall i \in I$. The cardinality $ \# I $ of $I$ is called {\bf size of $f$}. 

In \cite{HBG}, the {\bf Debreu dimension $d_D(X)$ of $X$} is defined by

$$ \dim_{D} (X) := \inf \{ \# \{ A \subset P( X \times X ) \} \vert \forall a \in A: a {\rm \ total \ Debreu \ separable \ order \ on \ } X {\rm  \ and \ } \leq = \bigcap A  \} ,$$ 

\V where an order on $X$ is called {\bf Debreu separable} iff there is a countable subset $N$ of $X$ such that $\forall x , y \in X: x \leq y \Rightarrow (\exists n \in N: x \leq n \leq y)$. Debreu \cite{gD} proved that a total order on $X$ is Debreu separable iff there is a utility of size $1$ on $X$. 

\bigskip

Furthermore, \cite{HBG} defines the {\bf Hack-Braun-Gottwald (HBG) dimension} $\dim_{{\rm HBG}}$ {\bf of $X$} as

$$ d_{{\rm HBG}} (X) = \inf \{ \# I \vert \exists f: X \rightarrow \R^I {\rm \ utility \ on \ } X \} ,$$ 

in other words, the HBG dimension is the minimal size of a utility on $X$. Ore's description of $d_{DM}$ shows that $d_{{\rm HBG}} \geq d_{DM}$ as the totally ordered spaces $C_i$ are restricted to be order equivalent to $\R$ in the HBG dimension.
Of course, by definition $\dim_g\geq \dim_D$, as the chains in the HBG dimension do not have to define total orders on $X$. That is, all in all we have $ \dim_{DM}  \leq \dim_g \leq \dim_D$. In \cite{HBG}, the lexicographical $\R^2$ is shown to have Dushnik-Miller dimension $1$ and uncountable HBG dimension.\footnote{Indeed, $X:= \R^n_{{\rm lex}}$ does not admit a utility $f$ of size $1$, as $f$ would map each $\{ x \} \times \R $ to a proper interval $I_x$, and the uncountably many $I_x$ would be pairwise disjoint. Now, if there was a countable-sized utility $\{ f_n \vert n \in \N \}$ on $X$, then we could define $F:= \sum_{i=1}^\infty 2^{-n} \arctan \ci f_i : X \rightarrow \R$ which is stricly increasing: For $x \leq y $ we have $f_i(x) \leq f_i(y)  \ \forall i \in \N$, thus if not $F(x) < F(y)$ we get $f_i(x) = f_i (y) \ \forall i \in \N $ implying $x=y$. And finally, each strictly increasing function on a totally ordered set is a utility.} and an example of an ordered set $X$ with $d_{{\rm HBG}} (X) = 2$ and $d_D(X) = \N$ is given. 

\bigskip

In any case, we never leave the reign of sets, as in the Dushnik-Miller and Debreu case, from the definition, $\# I \leq \# P(P(X))$, whereas in the case of HBG dimension w.l.o.g. $\# I \leq \# P(X)$.

The Dushnik-Miller dimension can be localized e.g. by infimizing over order intervals containing a given point $p$, then denoted by $\dim_{DM,l} (M,p)$. By definition, each localized dimension $\dim_{\cdot , l} $ is less or equal to $\dim_{\cdot}$. Still so, Dushnik-Miller dimension, HBG dimension and Debreu dimension, whether localized or not, all fail to recover the manifold dimension even in simple cases, as on Minkowski space of dimension $\geq 3$ all of them take at least the value $\aleph_0$: 

\begin{Theorem}
	
	$\forall n \setminus \{ 0,1\} : \dim_{DM} (\R^{1,n})  = \dim_{g} (\R^{1,n})  = \aleph_0$,
	
	$\forall n \setminus \{ 0,1\} \forall p,q \in \R^{1,n} : \dim_{DM} (J(p,q))  = \dim_{g} (J(p,q))  = \aleph_0,$
	
	$ \dim_{D} (\R^{1,n}) , \dim_{D} (J(p,q)) \geq \aleph_0$.
\end{Theorem}

\V{\bf Proof.} For $\dim_{DM} (\R^{1,n}) \geq \aleph_0$, we show that the dimension is not finite as then there would be $p \in \partial J^+(0) $ such that $(\partial J^+ (0)) \cap (\partial J^+ (p))  $ is not totally ordered. More precisely, for $X=\R^{1,n}$, the image of $X$ in all chains $C_i$ in a minimal Ore product have to be order equivalent: Each two points in $X$ can be joined via a piecewise causal curve $c$ (this is even true in any connected Lorentzian manifold). Now $f_i \ci c$ is piecewise monotonous and continuous in the interval topology on $X$ (which just coincides with the usual product topology) and thus on $C_i$. Thus every chain is connected (complete). It is furthermore separable: The image $f_i(S)$ is dense for each separator $S$ in $X$. Thus $f_i(S)$ is order-isometric to a real interval. Now assume the dimension $m$ is finite. Then $f$ is an injective continuous map from $\R^{n+1}$ to $\R^m$. Take a point $p \in f(\partial J^+(0) \setminus \{ 0\} )$. Now let $p \in J^+(0) \setminus \{0\}$ such that $ d:= \# \{ i \in \N_m^* \vert f_i(p) \neq 0 \}\geq n >1$, such a point exists due to invariance of the domain. So let $F$ be a hyperplane of dimension $d \in [n; m)$ containing an open neighborhood $U$ of $f(p)$. Then $D_U := f(\partial J^+(0) \cap U) =  f(U) \cap F  $. The totally ordered subsets in $ \R^m$ and thus in $F$ are just the rays from $0$. But $f(D_U) $ cannot be contained in the ray through $f(p)$, again by invariance of the domain. On the other hand, $\partial J^+(p) \subset \partial J^+ (0)$ by the push-up property. Thus $\partial J^+ (p) \cap \partial J^+ (0) \cap U$ is not totally ordered, contradiction. Then use the classical fact (only using the countable axiom of choice) that the cardinality of $\N$ is smaller than the cardinality of each other infinite set. For $\dim_g (\R^{1,n}) \leq \aleph_0$, we easily show that $\leq = \bigcap_{v \in \mathbb{S}^{n-1} \cap \Q^{n}} (e_0 + v )^{-1} (\leq_\R) $.
\hfill \qed

\bigskip

In particular, there is no order embedding of $\R^{1,n}$ into $\R^m_{prod}$ for any $n,m \in \N$. The $\R^m_{prod}$ (which satisfy $\R^m_{prod} \times \R^n_{prod} = \R^{m+n}_{prod}$) arise in a very elementary manner as the causal structures of Lorentzian length spaces as follows: Let us describe causal cylinders over metric spaces, which are g.h. Lorentzian pre-length spaces generalizing standard static spacetimes: Let $(M,d)$ be a metric space. Then the {\bf causal cylinder $C(M,d)$ over $M$} is defined as $X:= \R\times M$ with Lorentzian distance $\sigma_d : X \times X \rightarrow \R  $ given by 

$$ \sigma_d ((t,p), (s,q)) :=  \cdot \sqrt{(s-t)^2 - d^2(p,q)} \ {\rm \  for \ \ }  s - t  \geq d(p,q)$$  

and $- \infty $ otherwise. If $(M,d)$ is a length space then $C(M,d)$ is a Lorentzian length space. (This definition coincides with the notion of generalized cone as in \cite{sAmGmKcS} specialized to the case $f =1$.) For a real interval $I$ we denote the corresponding subset $I \times M$ of $C(M,d)$ by $C_I(M,d)$. Then we can identify $\R^n_{prod} = C ((\R^{n-1}, d_{\ell^\infty}))$, another way to show that $\R^n_{prod}$ is a Lorentzian length space.

\bigskip

In the category of Lorentzian length spaces, or of ordered measure spaces, there is the Lorentzian Hausdorff dimension as defined\footnote{During the preparation of this article, the author learned from Sumati Surya that this notion actually appears under the name of "midpoint scaling dimension" in works of Sorkin, cf. \cite{rS}, p.11, however without a precise statement about it.} in \cite{rMcS} or \cite{oM:LorFunct}. Here, however, we try to assign a dimension to ordered spaces {\em without a specified measure}, i.e., dimensions for the category $C_8$ of ordered sets. An underlying question is how large the space of isomorphism classes of $n$-dimensional globally hyperbolic (g.h.) spacetimes is (w.r.t. the Lorentzian Gromov-Hausdorff metrics) in the space of isomorphism classes of $n$-dimensional compact g.h. Lorentzian length spaces satisfying a bound on curvature and Lorentzian injectivity radius, which could provide insights on the question of Th. \ref{InversesOnClosure} to which extent the functors between the categories of g.h. compact ordered measure spaces and g.h. compact Lorentzian length spaces from \cite{rMcS} and \cite{oM:LorFunct} are mutually inverse on subsets with such bounds. In the rest of this section we assume $\ll = \beta (\leq)$ and $ X = \underline{X} := \{ p \in X | I^+(p) \neq \emptyset \neq I^-(p)\}$. For each presumptive functor $\dim_*$ defined below, we ask

\begin{enumerate}
	\item whether it is a dimension in the sense of Eq. \ref{Dim-Def};
	\item whether it is {\bf compatible with the manifold dimension}, i.e. whether for each $n$-dimensional g.h. Lorentzian spacetime $(M,g)$ and $X:= (M, \leq)$ we get $\dim_* (X,p) = n$ for all $p \in X$;
\end{enumerate}

\V{\bf Remark.} Obviously, no compatible dimension can be upper or lower semi-continuous w.r.t. GH topology: We easily find static counterexamples by collapse of the Cauchy surfaces, respectively by approximating a Riemannian Cauchy hypersurface by surfaces of increasing genus. Furthermore, we could not require the embeddings in \ref{Def-minimal} to be order embeddings, because, via the Malament theorem, the Weyl tensor provides a local obstruction.

\bigskip

A first approach to define dimension is taken from a preprint of Stoica \cite{cS}: With the usual horismotic relation $E^\pm := J^\pm \setminus I^\pm$, the {\bf horismoticity} $\de_h$ {\bf of $X$ at $p$} is defined by 

$$ \de_h (X,p) := \inf \{ \# A \vert A \subset X \land \bigcap_{a  \in A} E^\pm(a) = \{p \} \}.$$

\begin{Theorem}[Properties of horismoticity]
	\ 
	\begin{enumerate}
		\item $\de_h $ is no dimension in the sense of Eq. \ref{Dim-Def}.  
		\item $\de_h$ is compatible with the manifold dimension.
	\end{enumerate}
	
\end{Theorem}

\V{\bf Proof.} To the see the second assertion, for $q_i$ in the injectivity radius neighborhood $U$ of $p$, $f_i:= \s^2 (\cdot, q_i)$ is differentiable in $p$ (however not $C^1$) if $q_i \in \partial J^+(p)$, and then ${\rm grad}_p f_i = v_i:= c_i'(0)$ where $ c_i$ is the unique null geodesic in $U$ with $c(0)= p, c(1) = q_i$. Therefore for $n = \dim (M)$ points in general position, those vectors span the whole space, such that $p$ is the unique point on the intersection of all horismos. Now consider a standard static Lorentzian length space $X:= \R\times S$ (where $S$ is a complete length space). As the intersection of causal cones with standard slices are balls, we easily see that $X$ has horismoticity $n+1$ if and only if $\dim_B (S) = n$ where $\dim_B$ is {\bf Blumenthal's metric dimension of $S$} defined by   

$$ \dim_B (M,p) :=  \inf \{ \# A \vert A \subset M \land \forall q \in M : d(q,a) = d(p,a) \forall a \in A \Rightarrow p = q \} . $$

$\dim_B$ is compatible with the manifold dimension (an easy consequence of normal coordinates and convex neighborhoods).

 Regarding Item 1, \cite{sBaB} constructs simple examples of compact length spaces $A$ und $B$ (even metric graphs) with an injective isometric embedding $I: A \rightarrow B$ and $\dim_B(A) > \dim_B(B) $.  \hfill \qed

\bigskip

Another approach to define a notion of dimension is the {\bf minimal (compatible) dimension} $d_m$ defined as the minimal dimension with $\dim_m (\R^{1,n}) = n+1 \ \forall n \in \N$: For $N(p)$ being the neighborhood system of $p$, we define

\bea
\label{Def-minimal}
\dim_m (X,p) := \sup \{ \vert \exists V \in N_p \  \exists U \subset \R^{1,n} {\rm \ open \ } \exists f: U \rightarrow V {\rm \ injective \ increasing}    \}.
\eea

\begin{Theorem}[Properties of the minimal dimension]
	$\dim_m$ is monoidal. Moreover,
	\begin{enumerate}
		 \item $\dim_m$ is a dimension in the sense of Eq. \ref{Dim-Def}.
		 \item $\dim_m$ is compatible with the manifold dimension.
    \end{enumerate}	
	
\end{Theorem}

\V {\bf Proof.} To see that $\dim_m$ is a dimension, observe that given an monomorphism $\Phi: X \rightarrow Y$ of $C_8$, we can postcompose each monomorphism $f: \R^{1,n} \rightarrow X$ with $\Phi$ to get a monomorphism $\Phi \ci f : \R^{1,n} \rightarrow Y$. To see monoidality of the minimal dimension, look for a product cone fitting with its closure into the interior of the Euclidean cone at $0$ and extend a bit by continuity. 

To see compatibility, we observe that each monomorphism $M \rightarrow N$ between strongly causal spacetimes is continuous and then use invariance of the domain. On the other hand, let $x \in M$ be a point in a strongly causal spacetime, let $U$ be a geodesically convex and causally convex neighborhood of $x$. Via normal coordinates, we find an increasing map from a small neighborhood in the same-dimensional Minkowski space to $U$. \hfill \qed 

\bigskip
\bigskip

\V{\bf Remark.} The following definition of "catcher dimension" should be compared with the "catcher theorem" from \cite {oM:Hor}: Let $(M,g)$ be a spacetime with noncompact Cauchy surfaces. Then from any
point $p \in  M$ and any point $q \in I^+(p)$ there is a continuously inextendible timelike future curve $c$ from $p$ not intersecting $J^+(q)$. This is even true if we replace $q$ with a compact subset $Q\subset J^+(p)$, via the same proof.

\newpage

 Let $a \ll b \ll c \ll d$. We call $V \subset  J^+(a) \setminus I^-(c) $ a {\bf catcher$^\pm$ set for $(a,b,c,d)$} iff $J^\pm(b) \setminus I^\pm(V) \subset J^\mp(d)  $. Let $ C^\pm (a,b,c,d)$ be the set of catcher$^\pm$ sets for $(a,b,c,d)$, let $P$ be the set of partitions of $[0;1]$. For a timelike curve $c: [0;1] \rightarrow X$ and $p \in P$ define 
	
	\[ \dim^\pm _c (X,c,p) :=  \min \{ \overbrace{\max_{k \in \N_{\# p-2}} \# A_k}^{=:m(A)} | \overbrace{A_j \in C^\pm (p_j, p_{j+1}, p_{j+2}, p_{j+3}  ) \forall j \in \N_{\#p -2}    }^{\Leftrightarrow : A \in C^\pm(c,p) \subset P(X)^{\#p -2}} \} . \]

	$ \dim^\pm_c(X,c) := \liminf_{p \in P} \{ \dim^\pm_c (X,c,p) | p \in P \}$ and
	
	 $\dim_c^\pm (X,x) := \limsup_{(p,q) \rightarrow x} \sup \{ \dim (X,c) | c {\rm \ longest \ curve \ in \ } J(p,q)\}$.

	\begin{Theorem}[Properties of the catcher dimension]
		
		 $\dim^\pm_c \geq \dim_m$. Moreover,
		\begin{enumerate}
			\item  $\dim_c^\pm$ is a dimension in the sense of Eq. \ref{Dim-Def},
			\item $\dim_c^\pm$ is compatible with the manifold dimension.
		\end{enumerate}
	\end{Theorem}
	
	{\bf Proof.} The first assertion follows from Items 1 and 2. 
	
	For the first item, let $f:X \rightarrow Y$ be a monomorphism, let $x \in X$ and let $c$ be a temporal curve through $x$. Then $f \ci c$ is a temporal curve in $Y$ through $f(x)$, and for a tuple $A$ of subsets of $Y$ as in the definition of $\dim^\pm $, $f^{-1} (A)$ has the same cardinality and $f^{-1} (A) \in C^\pm (c, f^{-1}p)$.  
	
For compatibility: By considering small normal neighborhoods and the openness of the conditions, the statement is reduced to the case of Minkowski space. Let $x_0$ be the $0$-th coordinate function on Minkowski space and let $S_a := x_0^{-1} (\{ a\})$. We consider $J^+(0)$. For $p_i \in S_0 \setminus \{ 0\} \subset \R^{1,n} \setminus J^+(0)$, $J^+(p_i) \cap S_1$ is a disc of radius $ 1$. But indeed the $1$-ball around a standard $n$-simplex $S$ of diameter $1$ with midpoint $0$ contains the unit ball $B$. This follows from the fact that the union of half-spaces 
	
	$$ (p_i^\perp ) ^+ := \{ x \in \R^n : \langle x, p_i \rangle >0 \} = \{ x \in \R^n : d(0, \cdot) > d( p_i, 0)  \}$$
	
	is all of $\R^n$, as every point in the interior of $S$ is a positive linear combination of the vertices (in barycentric coordinates)  and a union of half-spaces contains with a vector $v$ also $\R^+ v$. This also holds for perturbations of $S$: there is a $\de >0$ such that the $1$-ball around each $\tilde{S}$ with $d_H(S, \tilde{S}) <\de$ contains $B$. Then, scaling invariance in Minkowski space allows to choose small simplices around each point of the partition. The other direction of the dimension estimate is as follows: Given $k <n$ points $q_1, ... q_k$ in $J(q^-,q^+) \setminus J(p)$, then for their orthogonal projections $v_1,..., v_k$ onto $e_0^\perp$, there is a vector $w$ at $p$ with $\langle w, v_j \rangle >0$ for all $j \in \N_k^*$. A small calculation shows that for equidistant points $p_1,...,p_4$ along a curve, each catcher set $A$ satisfies $\# A \geq n-1$. The same is true for {\bf $\e$-almost equidistant partitions} (defined by $d(p_i ,p_{i+1})/d(p_i, p_{i-1}) \in [1- \e, 1 + \e]$ for all $i$) for appropriate $\e >0$, and for each $\e>0$ and each partition can be completed to an $\e$-almost equidistant partition. 
	
	 \hfill \qed

\newpage

\section{Metrics on Cauchy sets and anti-Lipschitz Cauchy time functions}

We first establish a link between the convergence of causal cylinders and the convergence of their bases: 

\begin{Theorem}
	\label{FromCentersToCylinders}
	Let $(M_n , d_n) \rightarrow^{GH}_{n \rightarrow \infty} (M_\infty , d_{\infty})$ be a convergent sequence of compact metric spaces. Let $I \subset \R$ be a bounded interval and let $X_n := C_I (M_n)$. Then $X_n \rightarrow_{n \rightarrow \infty} X_\infty$ in $d^-_{GH}$ and $d^\times_{GH}$.  
\end{Theorem} 

\V{\bf Proof.} Let $n \mapsto \big( C_n : M_n \rightarrow M_\infty \big)$ be a sequence of correspondences with $\dist (C_n) \rightarrow_{n \rightarrow \infty} 0$. Then for $\widehat{C}_n := \mathrm{Id}_\R \times C_n = \{ ((r,m), (r,s)) \vert (m,s) \in C\}$ we calculate, for large $n$, 

$$  \dist (\widehat{C}_n) = \sup \{ \vert \sqrt{(t-s)^2 - d_n (p,q) } - \sqrt{(t-s)^2 - d_\infty (p,q) } \vert : (p,q) \in C_n \} \rightarrow_{n \rightarrow \infty} 0  $$

(as the case $  s - t  \leq d_{\infty} (p,q)$ has a zero contribution and the condition $  s - t  > d_{\infty} (p,q)$ is open). This concludes the proof for $d^-_{GH}$. For $d^{\times}_{GH}$ we first enlarge $I$ by $I_\e := (-\e; \e )$ show $d^-_{GH}$-convergence of $I_\e \times M_n$ and then apply Theorem \ref{ContinuousFunctor}. \hfill \qed

\bigskip

The (local at $p$) catcher dimension on a cylinder $\R\times C$ satisfies

\bea
\label{LebesgueToCompactness}
\dim_c (\R\times C) = \dim_M (C)  +1
\eea

\V if Minkowski's (local at $x:=  \pr_2(p)$) covering dimension $\dim_M(C,x)$ of $C$ exists, defined as 

$$ \dim_M (C,x) := \lim_{\de \rightarrow 0} \log_{1/\de} (N(C,\de))  $$

where $N(C, \de) $ is the minimal number of $\de$-balls covering $C$, resp. its local analogon. Equation \ref{LebesgueToCompactness} is easy to see: For an open neighborhood $U$ of $p$ we get that  $U_0 := U \cap (\{ \pr_1 (p)\} \times C)$ is an open neighborhood of $\pr_2 (p) $. If we can cover every $U_0 \setminus \{ p \}$ by $n$ balls of radius $\de$, the catcher dimension of $ [- \de; \de] \times C$ is less or equal $n$, and $\dim_{H} (\R\times C) = \dim_H(C) +1 $. 

\bigskip

From this and Th. \ref{FromCentersToCylinders} we conclude that the dimension of a Gromov-Hausdorff limit of $n$-dimensional spacetimes can be strictly less than $n$, even if all Cauchy surfaces have curvature uniformly bounded below --- just take the products of a collapsing metric sequence with a Lorentzian interval. Without a curvature condition, the dimension of a Gromov-Hausdorff limit of $n$-dimensional spacetimes can even be be strictly larger than $n$: We can GH-approximate any compact manifold by compact surfaces of increasing genus, and then we take the sequence of the products with the unit interval. Thus, this kind of collapse effects can occur in Lorentzian Gromov-Hausdorff convergence.

\bigskip

Fixing $p \in [1; \infty) $, we define a metric on each temporally compact Lorentzian length space $X$ by

$$  \forall x,y \in X:  \check D_p (x,y) :=  \vert \s_x^p - \s_y^p \vert_\infty  , 
$$

\V Let a {\bf Cauchy slab} be $J^+(S_-) \cap J^-_(S_+)$ for two connected spacelike Cauchy surfaces $S_-$ and $ S_+ \in I^+(S_-)$ of a g.h. spacetime. We find (for geodesicness adapt Prop. 2.5.19 and 2.5.22 from \cite{BBI})

\begin{Theorem}
	\label{NoldusP}
	Let $X$ is a Cauchy slab, $p \geq 2$ then $(X, D_p := \la (\check D_p))$ is a connected geodesic length space. \hfill \qed
\end{Theorem}

A small calculation shows that the intrinsification $D_p$ of $\check D_p$ on the strip $[-a;a] \times \R^n \subset \R^{1,n}$ equals the Euclidean metric w.r.t. the given coordinates. We denote $\check D_{p}$ as the {\bf Noldus$^p$-metric} and $D_p$ as the {\bf Noldus$^p$-length metric}.

\begin{Definition} Let $t$ be a real function on a Lorentzian length space $X$ (whose causal relation is denoted by $J$) carrying a metric $D$.
	\begin{enumerate}
		\item $t$ is called {\bf rushing} iff $t(y) - t(x) \geq \s (x, y)  $ for all $(x,y) \in J$.
		\item $t$ is called {\bf anti-Lipschitz (w.r.t. $D$)} iff $ t(x) - t(y) \geq D(x,y) $ for all $ (x,y) \in J$. 
		\item $t$ is called {\bf generalized Cauchy} iff there are $A_-, A_+ \in \R \cup \{ \pm \infty\}$ such that for every causal curve $c: (a_-;a_+) \rightarrow X$ not continuously extendible to a larger open interval we have $t(c(s)) \rightarrow_{s \rightarrow a_\pm} A_\pm$.
	\end{enumerate}
\end{Definition}

{\bf Remark.} The notion of "rushing" has been defined by E. Minguzzi in \cite{eM}. Obviously, any rushing function is a time function. As $\s (x,y) \leq D_{p} (x,y)$ for all $(x,y) \in J $ by the conditional inverse triangle inequality, any time function that is anti-Lipschitz w.r.t. $D_{p}$ is rushing. Obvious examples for $D_{p}$-anti-Lipschitz time functions on open subsets are distances $\Delta_S$ (w.r.t. $D_{p}$) from a given Cauchy set $S$. The interest in anti-Lipschitz Cauchy time functions also stems from the Sormani-Vega null distance whose definiteness is ensured if the used time function is anti-Lipschitz w.r.t. a local metric (see Lemma 4.3 in \cite{SV} for spacetimes and Prop. 3.8 in \cite{mKrS} for Lorentzian length spaces). For the following theorem we define, for $A \subset X$, a metric $D_p^A : X \times X \rightarrow [0; \infty]$ by $D_p^A (x,y) := \sup \{ | \s^p (x,r) - \s^p (y,r) | r \in A \}$. Obviously, $D_p^A < \infty$ for any temporally compact $A \subset X$.

\begin{Theorem}
	\label{RCF}
	Let $X$ be a Lorentzian length space with a Cauchy subset $S$, let $A \subset X $ such that $D:= D_p^A < \infty$. Then there is a $D$-anti-Lipschitz generalized Cauchy time function $T$ on $X$ with $S= T^{-1} (\{ 0\})$. 
	
\end{Theorem}

{\bf Proof.} We perform the proof only for the case that $I^+(x) \neq \emptyset \neq I^-(x) $ for all $x \in X$ and construct a Cauchy function; by an obvious generalization of the construction one obtains a anti-Lipschitz generalized Cauchy function in the case that there is a past and/or future boundary. The proof is done similar to the proof of existence of a steep Cauchy function as in \cite{oMmS}. We use the following building blocks:

\begin{enumerate}
	\item A Cauchy time function $t$ on $X$ whose existence is ensured by the result in \cite{aBlG},
	\item a {\em fat cone covering}: a well-behaved covering of each Cauchy surface by forward cones defined below,
	\item for each $ (p,q) \in J$ a continuous function $ \tau_{pq}: X \rightarrow [0; \infty)$ supported on each $J^+(p)$ and anti-Lipschitz on $J^+(q)$ defined below,
	\item the facts that sums of anti-Lipschitz functions and time functions are anti-Lipschitz and that for any continuous function $f$ on $X$ and any anti-Lipschitz function $s$ and any $K\subset X$ compact there is $C>0$ with $f + C \cdot s $ anti-Lipschitz on $K$ (proven as in \cite{oMmS}). 
\end{enumerate}

Let $S$, $S'$ be two Cauchy surfaces of $X$ with $S \in I^+(S')$. A {\bf fat cone covering of $X$ above $S'$} is a set of tuples of points $\{ (p_i,q_i) \vert i \in I \} \subset I^+ (S') \times I^+ (S')$ such that $ p_i \ll q_i \forall i \in I$ such that  $\{ I^+ (q_i ) \cap S \vert i \in I \}$ resp. $ \{ J^+ (p_i) \cap S \vert i \in I \}$ is a covering resp. a locally finite covering of $S$. The existence of a fat cone covering is proven in \cite{oMmS} using only facts available in the synthetic setting as well, essentially future compactness of $J^-(S)$ and continuity of $J^+$ and $I^+$. 

For fixed $(p,q) \in J$, we define $\tau_{pq} (x) := u (\s (p,x)) \cdot \sup \{ D(y,x) \vert y \in J(p,x) \} $ for a continuous increasing $u : [0; \infty ) \rightarrow [0; \infty]$ with $u (0) =0$ and $ u ([\s(p,q) ; \infty)) = \{ 1/\s (p,q)\}$, then it is easy to see that $\tau_{pq}$ satisfies all the requirements of Item 3 above, recalling $\s(p,x) \geq \s(p,q) \ \forall x \in J^+(q)$. 

Now, for two Cauchy subsets $Y,Z$ with $Y \leq Z$, given a fat cone covering of $Z$ above $Y$, denoting $\tau_i:= \tau_{p_iq_i}$, we can construct $\check \theta^+_{Y,Z}:= \sum_i c_i \tau_i : X \rightarrow \R^+ $ for appropriate $c_i \in (0; \infty)$ with support in $I^+(Y)$, time function on its support and anti-Lipschitz on $J^+(Z)$ with $\check \theta^+_{Y,Z} >1 $ on $J^+(Z)$. Analogously we define by time-duality $\check\theta^- (Y,Z) $, time function on its support contained in $I^-(Z)$ and anti-Lipschitz on $J^-(Y)$ with $\check \theta^-_{Y,Z} ,-1 $ on $J^-(Y)$. Now we define $\theta^+_{Y,Z} := \rho \ci \check\theta_{S',S}  $ for $\rho: \R\rightarrow \R$ continuous increasing with $\rho =0$ on $(- \infty; 0]$ and $\rho = 1 $ on $[1; \infty)$. Then $\theta^+_{Y,Z}$ is a time function on $(\theta^{+}_{Y,Z})^{-1} ((0;1))$ with $\theta^+_{Y,Z} =0 $ on $ J^-(Y)$ and $\theta^+_{Y,Z} =1 $ on $ J^+(Z)$. Next we choose Cauchy subsets $S^{--}, S^-, S^+, S^{++}$ with $S^{--} \leq S^- \leq S \leq S^+ \leq S^{++}$ and $I^+(p) \neq \emptyset \neq I^-(q)$ for all $p \in S^{++} $ and all $q \in S^{--}$ and put $\theta^- := \theta_{S^--, S^-}, \theta^+ := \theta_{S^=, S^++}$ and

$$ \theta := 2 \frac{(\Delta_S +1) \theta^-}{(\Delta_S+1)- \theta^+ +1} -1 . $$

Then $\theta $ is a time function on $J(S^{--}, S^{++})$, anti-Lipschitz on $J(S^-,S^+)$ with $\theta = -1$ on $J^-(S^{--})$ and $\theta = 1 $ on $J^+(S^{++})$. Our final definition is

$$ T:= \check\theta^-_{S^-,S} + \theta + \check\theta^+_{S, S^+}  $$

which eventually satisfies all the requirements of the theorem.  \hfill \qed

\bigskip

Of course, there is in general no continuous projection on the level sets of a rushing time function, nor is there, in general, some related product decomposition.

\bigskip

Let $(X,d)$ be a metric space and $S \subset X$. Then the "intrinsification of $d$ along $S$" $\la_S(d)$ is the induced extended length metric on $S$ (admitting the value $\infty$). We need semi-continuity of $\la_S$:

\begin{Theorem}
	\label{LengthizationSemicont}	
	$\la_S $ is lower semi-continuous w.r.t. the Gromov-Hausdorff metric on both sides, but in general not upper semi-continuous. 
\end{Theorem}

{\bf Proof.} Let us first show pointwise convergence, i.e. we fix two points $x,y$. Let $X_n \rightarrow_{n \rightarrow \infty} X_\infty$ in the Gromov-Hausdorff metric and let $\ell(c) \in ( \la (d_\infty) (x,y) - \e ; \la (d_\infty ) (x,y)) $. Then there is a polygonal arc $P$ with break points $P_1,...,P_N$ with  

$$ \sum_{k=1}^{N-1} d_\infty (P_k, p_{k+1}) = \ell(P) \in (\la (d_\infty (x,y ) - 2 \e); \la (d_\infty (x,y))) .$$

Furthermore, let $n \in \N$ such that there is a correspondence $C: X_n \rightarrow X_\infty$ with $ \dist (C) < \e/N$ and let $Q_k \in C^{-1}(P_k)$, $X \in C^{-1} (x)$, $Y \in C^{-1} (y)$. Then $\vert d_\infty(P_k, P_{k+1}) - d(Q_k,Q_{k+1}) \vert < \e/N$ for all $k \in \N_N$ and thus $d_n(X,Y) > d_\infty (x,y) - 3 \e$.
Compactness ensures that this is uniform: Let $N$ be a finite $\e$-net in $(X_\infty, \la (d_\infty)) $ of cardinality $N$, say, then we only need to perform the previous step $ N(N-1)$ times to show upper semi-continuity.

That $\la_S$ is in general not upper semi-continuous can be seen by applying $\la_S$ to the quotient metrics of $g_p (x,y) := |x - y |^p$ on $\mathbb{S}^1 = \R / \Z $ obtaining the discrete extended metric (taking the values $0$ and $\infty$) for $p \in (0;1) $ whereas $g_p \rightarrow_{p \rightarrow 1} g_1$, the standard metric. \hfill \qed

\bigskip

For a Cauchy set $S$ in a spacetime (considered as an almost LLS), the Noldus$^p$ metric on $S$ is in general very different from the induced Riemannian metric. There is a metric doing better in this respect: Let $S$ be a Cauchy set in an LLS $X$, let $D_r:= \{ p \in X : \max \{ | \tau(p, s) | : s \in S \} < r \}$. 
For a curve $c: I := [a;b] \rightarrow  S$, a partition $ P:= (t_0 = a, x_1, ..., x_n = b) $ of $I$ and $ p:= (p^+ , p^-) \in (D_r^n)^2  $ we define $p \in CD_r (c,P) : \Leftrightarrow  c(t_m) \in J(p_{m-1}^- , p_{m-1}^+) \cap J(p_m^- , p_m^+) \forall k \in \N_n^*$ (we call $p \in CD_r(c,P)$ a {\bf chain of diamonds for $(c,P)$ at distance $r$}). Let $\ell_r (c,P) := \inf \{ \sum \s (p_i^- , p_i^+) | p \in CD_r(c) \}$ (decreasing in $r$), $\ell (c) := \sup \{ \ell_r (c,P)| r>0, P {\rm \ partition \ of \ } I \} \in \R\cup \{ \infty \} $ and $d_S(x,y ) := \inf \{ \ell (c) | c: x \leadsto y \} \forall x,y \in S$.

\begin{Theorem}
	\label{FactsOnDelta}
	
	\ 
	\begin{enumerate}
		\item $(S, d_S)$ is a pseudometric space, whose induced topology is weaker than the relative topology of $S$.
		\item If $X$ is a spacetime and $S$ is a differentiable spacelike hypersurface, then $d_S$ coincides with the metric induced by the induced Riemannian metric on $S$.
		\item $d_{S} \leq 2 \la_S (\check D_1)$ (in particular, if $\la_S(\check D_1) < \infty$, then $d_{S} < \infty$).
		\item If $S_n \rightarrow_{n \rightarrow \infty} S_{\infty}$ in the pointed $\s$-Hausdorff metric, then $(S_n, d_{S_n}) \rightarrow (S_\infty, d_{S_\infty})$ in the pointed Gromov-Hausdorff sense.
	\end{enumerate}
\end{Theorem}

\V {\bf Remark.} In general, $d_S$ is not finite: Consider the example of an LLS induced by a discontinous Lorentzian metric conformally equivalent to Minkowski space whose conformal factor has a sufficiently integrable pole, more precisely, let $f: \R \setminus \{ 0\} \rightarrow \R$ with $f(s) := 1/|s|$, and $V_n := (\R^{n+1}, (f \ci \pr_0) \cdot g_{1,n})$. Then for each $x \in \R^n$ we have $\s ((0,x), (r,x)) = 2 \sqrt{x}$ and for each $s,t>0$, for $c_s (t) := (0,s+t, 0 ,..., 0) $ we get $\ell_1 (c) :=  2 \sqrt{t}$, and the intrinsified length of the segment between $ 0$ and $(0,1,0...,0)$ is infinite. The same construction with $f$ replaced with $g: s \mapsto s^2$ shows that $d_S$ is not a metric in general (In this case $d_S$ is identically $0$ on $\{ x_0 = 0\}$).

\V{\bf Proof.} For Item (1) we first show that the subset $M_X$ resp. $M_S$ of those points of $X$ resp. of $S$ that are in the image of the interior of the domain of definition of a longest curve is dense in $X$ resp. in $S$ (this is obvious by the fact that $\tau_+$ on the subset of points with $I^\pm (x) \neq \emptyset$ is the Alexandrov topology whose basis are causal diamonds, and by geodesic connectedness of globally LLLS). Now let $x \in S$ and $\e>0$. Then because of the above denseness of points on longest curves, there is such a point $q \in B(x, \e)$ and a longest curve $c: (- \beta; \beta) \rightarrow X$ with $q = c(0)$. For $\de := d(x,q)$ we choose $\mu \in (0; \max \{ \beta, \e - \de \})$, and then $J(c(- \mu ); c(\mu)) \cap S \subset B(x, \e)$.

For Item (2) we note that a length space metric $\la (d)$ is uniquely determined by the datum of $d$ in arbitrarily small neighborhoods. We consider Fermi coordinates around $S$, in which the metric is $-dt^2 + g_t$, and use local boundedness of $\dot{g}$. 

For Item (3), using denseness of points on longest curves as in Item 1, let $y \in M_S$ and choose $ q_0 \in \partial J^+(x) \cap I^+(y)$ such that the maximal extension of the longest curve $c_+: y \leadsto q_0$ intersects $\partial J^-(x)$ in a point $p_0$ (this is possible for $x$ close to $y$ by continuity of $J^\pm$). Then $J(p_0, q_0)$ contains $x,y$, and 

$$ \s (p_0, q_0)  = \ell (c_+ ) + \ell (c_-) = \s (q_0, y) + \s (y, p_0) = |\s(q_0,y) - \s (q_0,x)| + |\s (p_0, y) - \s (p_0, x)| \leq 2 \ccD(x,y),$$ 

thus $d_S (p_0, q_0) \leq  \ccD_1 (x,y) $ which concludes the proof by continuity and denseness. \hfill \qed

\bigskip 

Note that intrinsification preserves separability and local compactness if we additionally assume curvature bounds. An interesting
open question in this context is whether there is a good analogon for the {\em second} fundamental form of Cauchy surfaces in the synthetic setting.

\medskip

{\bf Remark.} There is a different approach to induce distances on Cauchy hypersurfaces of spacetimes in a synthetic manner using volumes \cite{EVS}, but the calculations there related to the recovering the true Riemannian metric seem to be probabilistic in nature, assuming a certain stochastical process.

\begin{Theorem}
	Let $p \in [1;\infty )$, let $X$ be an almost Lorentzian length space and let $t$ be a $\check D_{p}$-anti-Lipschitz generalized Cauchy time function on $X$ with $I:= t(X)$. Then the map $s \mapsto (t^{-1} (\{ s\}), d_s := \de_{t^{-1}(s),p})$ is continuous when we apply the Hausdorff distance w.r.t. $\check D_p$ (i) and also for the Gromov-Hausdorff metric w.r.t. the restrictions  of $\check D_p$ (ii). It is lower semi-continuous if we consider the Hausdorff metric w.r.t. $D_p = \la_X (\check D_p)$ (iii) on the right-hand side, or the Gromov-Hausdorff metric $D_p = \la_X(\check D_p)$ (iv), or the one-parameter family of metrics $\{ \la_{S_a}(\check D_p) | a \in I \} $ (v). 
\end{Theorem}

\V{\bf Proof.} (i) is obvious by continuity of $\s$ and compactness of the Cauchy surfaces, and implies (ii). Item (i) implies with Th. \ref{LengthizationSemicont} item (iii), which in turn implies (iv). Item (ii) implies with Th. \ref{LengthizationSemicont} Item 5. \hfill \qed

\bigskip

If $c$ is an acausal curve then the $\de_{p}$-length of the image of $c$ is the Lorentzian analogue of the $p$-packing measure, and it is well-known (cf \cite{Falconer}) that the packing dimension is greater or equal to the Hausdorff dimension.   
Furthermore, we have a nice dimensional relationship:

\begin{Theorem}
	Let an almost Lorentzian length space $(X, \sigma)$ have constant Lorentzian Hausdorff dimension $n$. Let $S$ be a Cauchy set of $X$. Then
	\begin{enumerate}
		\item $\dim_H (S,d_{S,p}) \leq n $.
		\item $S$ can be approximated in the Hausdorff metric w.r.t. $D_{p}$ by Cauchy sets of Hausdorff dimension $\geq n -1$. 
		\item If there is $K>0$ with $ \check D_{p} (x,y) \leq K \cdot \s (x,y) $ for all $(x,y) \in J$, then $\dim_{H} (S, d_S) \leq n-1$. In particular, this is satisfied if $X$ is the limit of Cauchy slabs with uniform bounds on timelike sec and timelike diameter.  
	\end{enumerate}
\end{Theorem}

\V{\bf Proof.} For $d \in [0; \infty)$, let $\a(d):= \frac{\pi^{d/2}}{\Gamma (d/2+1)}$ (being the volume of the $d$-dimensional unit ball for $d \in \N^*$) and $\omega(d) :=   \frac{\pi^{\frac{d-1}{2}}}{d \cdot \Gamma (\frac{d+1}{2}) \cdot 2^{d-1}}$. Recall 

$$\la_N^+ (U) := \alpha (N) \cdot \sum_i (\frac{\diam(U_i)}{2})^N $$ 

\V for $U \in CC^+ (A)$ and $\la_N^- (U) := \omega (N) \cdot \sum_i (\s (p_i,q_i))^N $ for $U \in CC^+ (A)$ (where $CC_\de (A)$ is the set of open coverings of $A$ with open sets of diameter $\leq \e$, for $CC_\de^-$ we additionally require each set to be a causal diamond $J(p_i, q_i)$).

On the other hand, we have 

$$\mu^\pm_{N, \de} (A) := \inf \{ \la_N^+ (U) | U \in CC^\pm_\de (A) \} $$

\V and $\mu_N^\pm (A) := \lim_{\de \rightarrow 0} \mu^\pm_{N, \de} (A)$. Then $\dim_H(A) = \sup \{ r \in [0; \infty) | \mu_r (A) = \infty \} = \sup \{ r \in [0; \infty) | \mu_r (A) = \infty \}$. Now each $U \in CC^-_\de (A)$ induces $U \cap \Sigma \in CC^+_\de (A \cap \Sigma)$, more precisely, for each $p_\pm  \in J^\pm (\Sigma) $ we get $\diam_{\check D_2} (J(p_-, p_+)) \geq \s (p_- , p_+) \geq \diam_d (J(p_-, p_+) \cap \Sigma)$.  This settles the first claim.

For the second assertion, let $(p,q) (n) \in CC_{1/n} (A)$ and $P(s,n) := \# \{ (p_k^{(n)}, q_k^{(n)} ) | t(p_k ) < s < t(q_k) \} / \# (p,q)(n) $. For generic $s$ there are $c,C >0$ with $cn \leq P(t,n) \leq Cn$. 

We consider a $D_p$-anti-Lipschitz Cauchy function $T$ on $X$ with $S = T^{-1} (\{ 0\})$ whose existence has been shown in Th. \ref{RCF} and $A := T^{-1} (- \e; \e)$. Then we apply Eilenberg's coarea theorem (\cite{BZ}, Ch.5) to the $1$-Lipschitz function $\de$ obtaining that 

$$  \mu_n (A) = \int_I  \mu_{n-1} (t^{-1} (\{ s\})) ds.$$

That is, $ \mu_{n-1} (t^{-1} (\{ s\}))$ is finite for almost all $s$. Furthermore, every such $(-\e; +\e)$ contains $s$ s.t. $T^{-1}(s)$ has open subsets of nonzero (but finite) $(n-1)$-dimensional $D_p$-Hausdorff volume, which consequently are of finite $(n-1)$-dimensional $D_p$-Hausdorff volume as well due to Th. \ref{FactsOnDelta}, Item 3 (however, the $(n-1)$-dimensional $\de_p$-Hausdorff volume of open subsets of $T^{-1} (s)$ may be zero as long as we do not have an estimate between $\de_p$ and $D_p$ in the other direction, as in the last item of the present theorem). This settles the second claim.

For the last assertion we only have to show that $\mu_r (S) = 0 $ for all $r > n-1$. Then we take $U_n \in CC_{1/n}^- (X) $ and observe that by the supposed estimates the part $W_n$ of the sum $\Sigma_n $ corresponding to subsets intersecting $S$ nontrivially is in $O(1/n) \cdot \Sigma_n$, thus for $r > n-1$ we have a constant $K>0$ with $ W \leq K \cdot 1/n \cdot \Sigma_n \leq \alpha (N) \cdot \sum_i (\frac{\diam(U_i)}{2})^{r+1} $, which converges to $0$. \hfill \qed

\bigskip

Our last question is whether we can find a good replacement for the projection onto Cauchy hypersurfaces and the product decomposition for g.h. spacetimes. Let $S$ be a subset of a LLS $X$. Then we call a timelike curve $c: [0;1] \rightarrow X$ starting at $S$ {\bf orthogonal to $S$}, in symbols $c \perp S$, iff $\s (c(0), c(1)) = \s (S, c(1))$. We can define a relation $r_{ut} : S_u \rightarrow S_t$ by the orthogonal curves on $S_u$ which is obviously right-total. To the author, it is not clear whether it is also left-total, i.o.w., whether at every point of $S_u$ an orthogonal curve starts. The function $x \mapsto \sqrt{|x|}$ provides an example of a Hölder continuous hypersurface $S$ of $\R^2$ containing a cusp $(0,0)$ that is not initial point of any $S$-minimal geodesic to one of its sides. But at least for Cauchy hypersurfaces in spacetimes with a metric of regularity $C^{1,1}$ and $|u-t|$ small, left-totality is satisfied:

\begin{Theorem}
	Let $P$ be a Cauchy surface in a globally hyperbolic spacetime with a metric of regularity $C^{1,1}$. Then each point $q \in P$ is initial point of a locally $P$-maximizing future and past geodesic.
\end{Theorem}

\V{\bf Proof.} Every semi-Riemannian submanifold $P \subset M$ has a normal neighborhood $V$ in $M$, i.e., there is an open neighborhood $U$ of the zero section $Z$ over $P$ s.t. for the normal exponential map $\exp^\perp$ we have $E:= \exp^\perp \vert_U$ is a diffeomorphism onto $V$ (see e.g. \cite{oN}, Prop.7.26). In particular, $E^{-1} (P) = Z$, and in local coordinates $\{ x_i \vert i \in \N_n \}$ around $q \in P$, the curve $ t \mapsto (t, x_1(q),..., x_n(q) )$ is geodesic, and $g_{0i} (q) = 0 \ \forall i \in \N_n^*$ (as $d_0 \exp^\perp $ is the canonical identification $T_0T_qM \rightarrow T_qM$). We will need regularity at least $C^1$ in order to make sure that geodesics are defined. Like in the smooth case where we would argue with the first variational formula, for general $C^1$ metrics we still get that each $P$-maximal geodesic is initially orthogonal to $P$ (\cite{mG}, Cor.2.5). For $C^{1,1}$ metrics, the connection is of Lipschitz regularity which implies that geodesics do not split and that locally maximal curves are geodesics \cite{eM:Sprays}. This, together with the fact that the exponential map is a diffeomorphism, implies that every $q\in P$ is the initial point of a $P$-maximal curve.   \hfill \qed

\bigskip

Finally, there is an analogon to the local-in-time boundedness of $|\dot g_t|$ on spatially compact g.h. spacetimes:

\begin{Theorem}
	Let $t$ be a $D_p$-anti-Lipschitz Cauchy function on a spatially compact Lorentzian length space $X$. Then the map $ a \mapsto (t^{-1} (a), \check{d}_{t^{-1} (a)})$ is continuous in the Gromov-Hausdorff metric, and the map $ a \mapsto (t^{-1} (a), d_{t^{-1} (a)})$ is lower semi-continuous in the Gromov-Hausdorff metric. 
\end{Theorem}

\V{\bf Proof.} Let $\e >0$. By compactness of $S_a$, there is $\rho >0$ such that for every $x \in S_a$ there are $x^\pm \in J^\pm (x) $ with $\s (x^-, x^+) > \rho $. Let $M_a:= d^{-1} ([0, \rho]) \subset S_a \times S_a$. For all $(x,y) \in M_a$ we choose $p^\pm (x,y)$ with $x,y \in D(x,y) := I(p^- (x,y), p^+(x,y))$ and $\s (p^-(x,y) , p^+(x,y)) \geq \check{\de} (x,y) - \e$. The set $\{ J(x,y) | (x,y) \in M_a\} $ is an open covering of $M_a$ and has a finite subcovering $D:= \{ D_1,..., D_N \} $ due to compactness of $M_a$. Then by causal continuity there is a $\de >0$ such that $D$ is an open covering of $M_{a+\de}$. We define a correspondence between $S_a$ and $S_{a+\de}$ by being contained in a same $D_k$; this correspondence is easily seen to have $\de$-distorsion $< \e$. \hfill \qed

\bigskip

{\bf Acknowledgements.} The author would like to thank Tobias Beran and Jona Röhrig for useful comments on a previous version of this article.

\bigskip



\begin{thebibliography}{99}
	
	
	\bibitem{sAmGmKcS}
	Stephanie B. Alexander, Melanie Graf, Michael Kunzinger, Clemens Sämann: {\em Generalized cones as Lorentzian length spaces: Causality, curvature, and singularity theorems}. arXiv:1909.09575
	
\bibitem{sBaB}
Sheng Bau, Alan F. Beardon: {\em The Metric Dimension of Metric Spaces}	Comput. Methods Funct. Theory 13:295 --- 305 (2013)
	
	
	
	
	
	\bibitem{BBI}
	Dmitri Y. Burago, Yuri D. Burago, and Sergei Ivanov: {\em A course in metric geometry}, volume 33
	of Graduate Studies in Mathematics. American Mathematical Society, Providence, RI,
	2001.
	
	\bibitem{BZ}
	Yuri Dimitrievich Burago, Victor Abramovich Zalgaller: {\em Geometric Inequalities}. Springer-Verlag (1988)
	
	\bibitem{aBlG}
	Annegret Burtscher, Leonardo García-Heveling: {\em Time functions on Lorentzian length spaces}. arXiv:2108.02693
	
	
	\bibitem{gD}
	Gerard Debreu: {\em Representation of a preference ordering by a numerical function.} Decision processes, 3:159–165 (1954).
	
	\bibitem{bDeM}
	Ben Dushnik, E.W. Miller: {\em Partially ordered sets.} American Journal of Mathematics, 63: 600 --- 610 (1941)
	
	\bibitem{EVS}
	Astrid Eichhorn, Fleur Versteegen, Sumati Surya:
	{\em Spectral dimension on spatial hypersurfaces in causal set quantum gravity}, Class.Quant.Grav. 36, 23, 235013 (2019).  arXiv:1905.13498
	
	\bibitem{Falconer}
	Kenneth J. Falconer: {\em Fractal Geometry - Mathematical Foundations and Applications}.
	John Wiley (1997).
	
	
	
	\bibitem{nG1}
	Nicola Gigli: {\em Lecture notes on differential calculus on {\rm RCD} spaces}. arXiv: 1703.06829
	
	
	\bibitem{nG2}
	Nicola Gigli, Enrico Pasqualetto, Elefterios Soultanis: {\em Differential of metric valued Sobolev maps}. arXiv: 1807.10063
	
	\bibitem{mG}
	Melanie Graf: {\em Singularity theorems for $C^1$-Lorentzian metrics}. Commun. Math. Phys. 378 (2020), 1417-1450. arXiv:1910.13915
	
	
	
	\bibitem{HBG}
	Pedro Hack, Daniel A. Braun, Sebastian Gottwald: {\em On a geometrical notion of dimension for partially ordered sets}, arXiv:2203.16272v3
	
	\bibitem{dH}
	David A. Herron: {\em Gromov-Hausdorff distance for pointed metric spaces}. The Journal of  \\ Analysis 24, 1 --- 38 (2016) 
	
	
	\bibitem{mKcS}
	Michael Kunzinger, Clemens Sämann: {\em Lorentzian length spaces}.  Ann. Global Anal. Geom. 54, no. 3, 399 --- 447 (2018). arXiv: 1711.08990
	
	\bibitem{mKrS}
	Michael Kunzinger, Roland Steinbauer: {\em Null distance and convergence of Lorentzian length spaces}. arXiv:2106.05393
	
	\bibitem{aL}
	Alexander Lytchak: {\em On the geometry of subsets of positive reach.} manuscripta math. 115, 199–205 (2004)
	
	\bibitem{aLkN}
	Alexander Lytchak, Koichi Nagano: {\em Topological regularity of spaces with an upper curvature bound}, JEMS vol. 24, no 1, 137 --- 165. arXiv:1809.06183
	
	
	
	\bibitem{rMcS}
	Robert McCann, Clemens Sämann: {\em Lorentzian Hausdorff dimensions and measure}. arXiv:2110.04386
	
	\bibitem{eM:Sprays}
	Ettore Minguzzi:{\em Convex neighborhoods for Lipschitz connections and
		sprays}, Monatsh. Math. 177, 569-625 (2015). arXiv:1308.6675
	
	
	\bibitem{eM}
	E. Minguzzi. {\em Lorentzian causality theory}. Living Rev. Relativ., 22:3, 2019.
	https://doi.org/10.1007/s41114-019-0019-x.
	
	
	
	\bibitem{eMmS}
	Ettore Minguzzi, Miguel S\'anchez: {\em The causal hierarchy of spacetimes}. In H. Baum and D. Alekseevsky (eds.), vol. Recent developments in pseudo-Riemannian geometry, ESI Lect. Math. Phys., (Eur. Math. Soc. Publ. House, Zurich, 2008), p. 299 -- 358, ISBN=978-3-03719-051-7. arXiv:gr-qc/0609119
	
	
	
	
	
	
	\bibitem{eMsS}
	Ettore Minguzzi, Stefan Suhr: {\em Lorentzian metric spaces and their
		Gromov-Hausdorff convergence}. arXiv: 2209.14384
	
	\bibitem{MS}
	Andrea Mondino, Stefan Suhr:
	{\em An optimal transport formulation of the Einstein equations of general relativity}, J. Eur. Math. Soc. (2022). arXiv:1810.13309
	
	
	
	
	
	
	\bibitem{oM:Hor}
	Olaf M\"uller: {\em Horizons}. Advances in Theoretical and Mathematical Physics vol. 19 no 4 pp. 747---760 (2015). arXiv: 1111.4571
	
	\bibitem{oM:LorFin}
	Olaf Müller: {\em Lorentzian Gromov-Hausdorff theory and finiteness results}, General Relativity and Gravitation
	54:117 (2022). arXiv: 1912.00988
	
	
	\bibitem{oM:Compl}
	Olaf Müller: {\em Topologies on the future causal completion}, arXiv:1909.03797
	
	\bibitem{oM:LorFunct}
	Olaf Müller: {\em Functors in Lorentzian geometry: Three variations on a theme}. General Relativity and Gravitation volume 55, Article number 39 (2023). arXiv:2205.01617
	
	
	
	
	\bibitem{oMmS}
	Olaf Müller, Miguel S\'anchez: {\em Lorentzian manifolds isometrically embeddable in $L^N$}, Trans. Amer. Math. Soc. 363 (2011), 5367-5379.
	arXiv:0812.4439
	
	
	\bibitem{kN}
	Koichi Nagano: {\em A volume convergence theorem for Alexandrov spaces with curvature bounded above}, Mathematische Zeitschrift vol. 241, 127 --- 163 (2002)
	
	
	
	\bibitem{oN}
	Barrett O'Neill: {\em Semi-Riemannian Geometry With Applications to Relativity}, Elsevier (1983)
	
	
	
	
	
	\bibitem{Noldus}
	Johan Noldus: {\em The limit space of a Cauchy sequence of globally
		hyperbolic spacetimes}. arXiv:gr-qc/0308075
	
	
	\bibitem{oO}
	Oystein Ore: {\em Theory of graphs}, American Mathematical Soc. vol.38 (1987)
	
	\bibitem{rS}
	Rafael D. Sorkin: {\em Causal Sets: Discrete Gravity (Notes for the Valdivia Summer School)}, in: Proceedings of the Valdivia Summer School, edited by A.
	Gomberoff and D. Marolf (2003). arXiv: 0309009
	
	\bibitem{SV}
	Christina Sormani, Carlos Vega: {\em Null distance on a spacetime}. Classical and Quantum Gravity, 33, no 8 (2016). arXiv:1508.00531 
	
	
	\bibitem{cS}
	O.C.Stoica: {\em Spacetime causal structure and dimension from horismotic relation}. arXiv:1504.03265v2
	




\bibitem{Stefan}
Stefan Suhr: {\em Theory of optimal transport
for Lorentzian cost functions}. M \"unster J. of Math. 11, 13 --- 47 (2018).
arXiv:1601.04532




\end{thebibliography}
\end{document}